\documentclass[11pt]{amsart}
\usepackage{amsmath, amsthm, amsfonts, amssymb, graphicx}
\usepackage{natbib}
\allowdisplaybreaks

\numberwithin{equation}{section}
     \newtheorem{thm}{Theorem}[section]

\theoremstyle{definition}

\theoremstyle{remark}
     \newtheorem{rem}{Remark}[section]

\def\be{\begin{equation}}
\def\ee{\end{equation}}

\begin{document}

\baselineskip=18pt

\title[A local analysis  of  the axi-symmetric Navier-Stokes
flow]{A local analysis  of  the axi-symmetric Navier-Stokes
flow  near a saddle point and no-slip flat boundary 
}

\author{Pen-Yuan Hsu}
\address{Department of Mathematics, Tokyo Institute of Technology, Meguro-ku, Tokyo 152-8551, Japan}
\email{pyhsu@math.titech.ac.jp }
\author{Hirofumi Notsu}
\address{Waseda Institute for Advanced Study, Waseda University, 3-4-1, Ohkubo, Shinjuku-ku, Tokyo 169-8555, Japan}
\email{h.notsu@aoni.waseda.jp}
\author{Tsuyoshi Yoneda}
\address{Department of Mathematics, Tokyo Institute of Technology, Meguro-ku, Tokyo 152-8551, Japan}
\email{yoneda@math.titech.ac.jp}

\maketitle

\begin{abstract}
As one of the violent flow, tornadoes occur in many place of the world. In order to reduce human losses and material damage caused by tornadoes, there are many research methods. One of the effective methods is numerical simulations such as the work in a recent article \cite{IOT}.  The swirling structure is significant both in mathematical analysis and the numerical simulations of tornado.
In this paper, we try to clarify the swirling structure. More precisely, we do numerical computations on axi-symmetric Navier-Stokes flows with no-slip flat boundary. We compare a hyperbolic flow with swirl and one without swirl and observe that the following phenomenons occur only in the swirl case:
 The distance between the point providing the maximum velocity
magnitude $|v|$ and the $z$-axis is drastically changing around some time
(which we call it turning point). An ``increasing velocity phenomenon'' occurs near the boundary and the maximum value of $|v|$ is obtained near the axis of symmetry and the boundary when time is close to the turning point.
\end{abstract}

\section{Introduction}\label{s:intro}
In this paper we consider    a local behavior of the  3D-Navier-Stokes  flow near a saddle point (with hyperbolic flow configuration) and no-slip flat boundary. 
The Navier-Stokes  equations with no-slip flat boundary are expressed as 
\begin{eqnarray}\label{NS}
\partial_tv+(v\cdot \nabla)v-\nu\Delta v+\nabla p&=&0\quad \mbox{in}\quad \mathbb{R}^3_+\times[0,T),\\
\nonumber
v_0=v|_{t=0},\quad v|_{\partial\mathbb{R}^3_+}=0,\quad \nabla\cdot v&=&0\quad\mbox{in} \quad\mathbb{R}^3_+\times[0,T),
\end{eqnarray}
where $v$ is a vector field representing velocity of the fluid, and $p$ is the pressure. The term ``hyperbolic flow configuration'' which used here and after means there is $\delta>0$ (depending on $t$) such that
$v(t,x)\cdot e_z>0$, $v(t,x)\cdot e_r(x)<0$,
or $v(t,x)\cdot e_z<0$, $v(t,x)\cdot e_r(x)>0$  for $0<|x_h|<\delta$ and $0<x_3<\delta$,
where $e_z=(0,0,1), ~e_r(x)=(x_1/|x_h|, x_2/|x_h|,0)$ and  $|x_h|=\sqrt{x_1^2+x_2^2}$. 
For the definition of 2D-hyperbolic flow configuration, we just need a simple modification.

Recently, the Euler flow with a saddle point (with hyperbolic flow configuration) has been extensively studied.
  \cite{BL} obtained strong local ill-posedness results in the Sobolev spaces $W^{n/p+1,p}$ for any $1<p<\infty$ 
and in the Besov spaces $B^{n/p+1}_{p,q}$ with $1<p<\infty$ and $1<q\leq\infty$ and $n=2$ or $3$ by using a combination of Lagrangian and Eulerian techniques with a saddle point (with hyperbolic flow configuration) structure. 
 In particular, they settled the borderline Sobolev case $H^{n/2+1}$. 
After that \cite{EM} and \cite{BL1} produced 
similar results in $C^1$ case (also $C^m$ case).
On the other hand, \cite{KS}  (see also \cite{IMY, KZ, X} for related topics)
showed
 a two-dimensional Euler flow with a saddle point (with hyperbolic flow configuration) in a disk for
which the gradient of vorticity exhibits double exponential growth in time for all times.
 Their estimate is known to be sharp, namely, the double exponential growth is the fastest possible growth
rate.
These results show  that the saddle point (with hyperbolic flow configuration) bring (some kind of) unstabilizing effects.

Now let us look back the history of Navier-Stokes equations briefly. Modern regularity theory for the solutions to the Navier-Stokes equations began with the works of  \cite{L} and \cite{H}.
They showed
the existence of a weak solution $v : [0 , \infty ) \times \mathbb{R}^{3} \rightarrow \mathbb{R}^{3} $ lying in the class of
$L^{\infty}(0, \infty ;
L^2(\mathbb{R}^3)) \cap L^2(0, \infty ; \dot{H}^1(\mathbb{R}^3))$ which satisfies the global energy inequality,
where $\dot{H}^1$ is a homogeneous Sobolev space with degree one. 

After that,  \cite{P}, \cite{S},
\cite{La}, and their joint efforts lead to the following
 Prodi-Serrin-Ladyzhenskaya criterion for the  Leray-Hopf weak solutions.

 \begin{thm}[Prodi-Serrin-Ladyzhenskaya]\label{heatequationpertubation}
 Let $v \in L^{\infty}(0, \infty ;
L^{2}(\mathbb{R}^{3}) ) \cap L^{2}(0, \infty ;\dot{H}^{1}(\mathbb{R}^{3}))$ be a
Leray-Hopf weak solution to (\ref{NS}), which also
satisfies $v \in L^{p}(0,T;
L^{q}(\mathbb{R}^{3}))$, for some $p, q$ satisfying $\frac{2}{p} + \frac{3}{q} = 1,$
with $q > 3$. Then, the solution $v$ is smooth and unique on $(0,T]\times \mathbb{R}^{3}$.
\end{thm}

It is also
worthwhile to mention that the exceptional case of $v \in
L^{\infty}(L^{3})$ 
 was finally established in the  work by
\cite{ISS}.
After the appearance of the Prodi-Serrin-Ladyzhenskaya criterion,
many different regularity criteria of solutions to (\ref{NS}) was
established by many researchers working in the regularity theory of
(\ref{NS}). 
For example, a regularity criterion along streamlines (characteristic curves) was 
constructed (see \cite{CY}).
Besides these, other
important works such as \cite{GHM}, in which type I blow up was excluded for solutions to (\ref{NS}) under a regularity condition on the vorticity direction in the half space to the case of the no-slip boundary conditon (see also \cite{CF}, which is the pioneer work in this field), and
\cite{CSTY, KNSS} for
axisymmetric solutions to (\ref{NS}) (the equation to the axi-symmetric case can be rewritten to (\ref{axisymmetricNS})) in the whole space, are attracting a lot of
attentions. 
We now recall that a singularity of a Navier-Stokes solution $v$ at a time $T$ is called Type I if
\begin{equation}\label{typeI}
\sup_x|v(x,t)|\leq \frac{C}{\sqrt {T-t}}
\end{equation}
for some $C>0$.
If the singularity of $v$ does not satisfy the condition (\ref{typeI}), we say Type II singularity
 (see \cite{KNSS} for the type of singularity). 
The axi-symmetric Navier-Stokes equations is expressed as 
\begin{eqnarray}\label{axisymmetricNS}
\nonumber
\partial_t u_r+u_r\partial_ru_r+u_z\partial_zu_r-\frac{u_\theta^2}{r}+\partial_r p&=&\Delta u_r-\frac{u_r}{r^2},\\
\partial_tu_\theta+u_r\partial_ru_\theta+u_z\partial_zu_\theta+\frac{u_ru_\theta}{r}&=&\Delta u_\theta-\frac{u_\theta}{r^2},\\
\nonumber
\partial_tu_z+u_r\partial_ru_z+u_z\partial_zu_z+\partial_zp&=&\Delta u_z,\\
\nonumber
\frac{\partial_r(ru_r)}{r}+\partial_zu_z&=&0,
\end{eqnarray}
where $u_r=u_r(r,z,t)$, $u_\theta=u_\theta(r,z,t)$ $u_z=u_z(r,z,t)$ and $\Delta=\partial_r^2+(1/r)\partial_r+\partial_z^2$.
The vector valued function  $v:=u_re_r+u_\theta e_\theta+u_ze_z$ with $e_r:=\frac{1}{\sqrt{x_1^2 +x_2^2}}(x_1, x_2, 0)$, $e_\theta:=\frac{1}{\sqrt{x_1^2 +x_2^2}}(-x_2, x_1, 0)$ and  $e_z=(0,0,1)$  represents velocity of the fluid, and $p$ is the pressure.


It is known that axi-symmetric solutions with no swirl (namely, $u_\theta\equiv 0$) have to be regular. See \cite{KNSS, La2, UI}.
On the other hand, \cite{CKN} showed that the axi-symmetric Navier-Stokes solution
could only blow up on the axis. 

In \cite{CHKLSY}, they considered a 1D transport equation with an additional variable which comes from a square of the swirl component $u _\theta$ of the velocity field in the 3D axi-symmetric Euler flow.
They showed that the transport equation can exhibit finite-time blow-up from smooth initial data using a contradiction argument, and their conclusion in \cite{CHKLSY} is that in the 3D axi-symmetric Euler flow with the flat boundary, the best chance for possible singularity formation seems to be at the saddle point (at the axis of symmetry intersect the boundary) (see also \cite{LH}).
On  the other hand, in \cite{K}, he constructed a regularity theory. More precisely,  suitable weak solutions of the three-dimensional axisymmetric Navier-Stokes equaions in a half space with no swirl case are Holder continuous up to the boundary exept for the origin.
His result suggests that even in the Navier-Stokes flow case, the best chance for possible singularity formation
may be at the saddle point.

Although there are many fruitful results based on mathematical analysis as we recalled above, it is not easy to analyze locally such fluid mechanics to go a step further mathematically. Thus, it should be effective for us to attempt numerical approach.  

Besides, from the above mathematical literature, the saddle point seems to be a key place and the hyperbolic flow with swirl might be  a key structure. Actually, the swirl ratios are very significant in the researches of tornadoes, see \cite{WD,IOT,N,IL}. There are several methods for researching tornadoes. Although the studies on real tornadoes and simulated tornadoes in laboratory are also important, numerical simulations provide a safe and cost effective way to analyze the behavier of tornadoes. We concentrate on the numerical computation in the swirl case and the behavior near the axis of symmetry (saddle point) at the boundary. From a mathematical point of view, the saddle point bring some kind of unstabilzing effects. In the point of view from the studies of tornadoes, the behavior insides the core of tornado (the center region near the $z$-axis) is  significant and might be very different to it outsides the core. The behavior near the ground (lower boundary) is also very significant for preventing the damage caused by tornadoes. Our numerical result can be compared with the mathematical literature (especially regularity results, see \cite{CF,GHM}) and also with the two-celled vortex structure in the tornado researches (especially numerical simulations). For more references, refer to \cite{GHM} for regularity results, refer to \cite{LH} for numerical studies of the Navier-Stokes and Euler equations, and refer to \cite{N,IL} for the studies of tornado-like vortices.

In the next section,
we will show the numerical computation on the difference between the swirl case and the no swirl case. We also show that  the maximum value of $|v|$ occurs at the place near the axis of symmetry (saddle point) and the boundary when the time approaches to the critical turning point (t=0.35 and t=1.0). (See figure 1,2,...,7, especially figure 7)

\section{Setting of the initial data and numerical results for the axi-symmetric Navier-Stokes flow}

In this section we set  an initial  data of  the axi-symmetric Navier-Stokes flow with a saddle point (with hyperbolic flow configuration).
We will compare two flows: with swirl ($u_\theta\not\equiv 0$)  and no swirl ($u_\theta\equiv 0$) cases.
In our numerical computation, we use the following cylindrical domain $\Omega$:
\begin{equation}\label{Omega}
\Omega:=\left\{x=(x_1,x_2,z)\in\mathbb{R}^3: 
-a<z<4a,
 \sqrt{x_1^2+x_2^2}<1\right\}
\end{equation}
and impose no-slip boundary condition:
\begin{equation}
v=0\quad\mbox{on}\quad \partial\Omega.
\end{equation}
 We set  the initial data in the following manner. (See also figure 1 and remark \ref{rem2.2}.)
Let $\varphi (a,\epsilon ,\sigma)=(a^2 +\epsilon)^{\sigma}$, we set the initial velocity for the swirl case as follows:
\begin{eqnarray}\label{initial}
u_z&=&\varphi (r, {\epsilon }_1, -{\beta}_1)\varphi (z, {\epsilon}_2, -{\beta}_2),\\
\rho&=&\varphi(r, {\epsilon}_3, -{\beta}_3)\varphi (z, {\epsilon}_4, {\beta}_4),\\
u_r&=&\mathrm{sign}(z)\rho u_z,\\
u_{\theta}&=&\varphi(r, {\epsilon}_5, -{\beta}_5)\varphi (z, {\epsilon}_6, -{\beta}_6)\label{initialend},
\end{eqnarray}
 where ${\epsilon}_i$ and ${\beta}_i$ ($i=1,2,...,6$) are constants. In the following numerical calculation, we set all ${\epsilon}_i$ and ${\beta}_i$ equal to $1$.
As for the no swirl case, we only change $u_{\theta}$ to zero.

\begin{rem}
By our setting of the initial data, the initial velocity magnitude $|v|$ at $(x_1,x_2,z)=(0,0,0)$ (we call this point the center of the initial velocity) is larger than it in other places. Phenomenons observed in this paper are more clear when the center of the initial velocity is close to the lower boundary (but not on the boundary). This is the reason for our choice of computational domain $\Omega$ like \eqref{Omega} instead of a symmetric one. See also appendix for comparison between the behaviors of two different centers.
\end{rem}

\par
 Let $\tau >0$ be a time increment of the computation and $h>0$ a
 (representative) mesh size.
 We perform the computation by a stabilized Lagrange-Galerkin (finite
 element) scheme~\citep{N1, NT1, NT2, NT3} based on the idea
\begin{subequations}\label{scheme}
\begin{align}
\frac{1}{\tau} \Bigl\{ v^k(x) - v^{k-1}\bigl( x-v^{k-1}(x)\tau \bigr)
\Bigr\} - \frac{1}{Re}\Delta v^k + \nabla p^k = 0,
\label{scheme:1st}\\
 \nabla\cdot v^k - \delta_0 h^2\Delta p^k =0,
\label{scheme:2nd}
\end{align}
\end{subequations}
 for $k=1, 2, \cdots$, where $v (x, k \tau)$ is simply denoted by $v^k(x)$, $Re$ is the Reynolds number, and the
stabilization parameter is set as $\delta_0=1$.
 We note that under some conditions the scheme has mathematical
 convergence properties
 of order $O(\tau +h)$ for the velocity in $H^1(\Omega)$ and for the pressure
 in $L^2(\Omega)$ and
 of order $O(\tau +h^2)$ for the velocity in $L^2(\Omega)$.
 The maximum, minimum and average mesh sizes are $1.88\times 10^{-2}$, $1.48\times 10^{-3}$ and $8.95\times 10^{-3}$, respectively,
 where the mesh size around $z$-axis is smaller than that of other part.
 The time increment is set as $\tau=1.25\times 10^{-2}$.
 In the following we show only numerical results by using the mesh size
 and the time increment mentioned above,
 where numerical results with a coarser mesh size and a larger time
increment are qualitatively similar. See appendix for the discussion
of the dependency of the numerical results on the discretization
parameters $h$ and $\tau$.
\begin{rem}\label{rem2.2}
 The initial velocities for the swirl and no swirl cases do not satisfy
 divergence free condition and the no-slip boundary condition,
 the computational velocities after the first time step, however, satisfy
 both conditions numerically, where the former condition is
satisfied in the sense that the equation \eqref{scheme:2nd} holds. The structure of initial data is useful for observing the swirling flow. It is also reasonable for us to consider the initial data \eqref{initial}-\eqref{initialend} as a pre-stage, the velocities after the first step as a real initial velocity in numerical computation. Although this construction of initial data is useful for observing phenomenons in numerical approach before the construction of more smooth initail data, a more careful choice of the initial data is desirable for mathematical analysis.
 \end{rem}
\begin{rem}
The axial symmetry is not explicitly imposed in the three-dimensional
computation by the stabilized Lagrange-Galerkin scheme, while the
problem setting including the initial velocity has the symmetry.
We perform the computation in order to show qualitative properties of
the effect of the swirl.
\end{rem}
\begin{figure}\label{fig1}
\includegraphics[width=7.0cm
,keepaspectratio
]{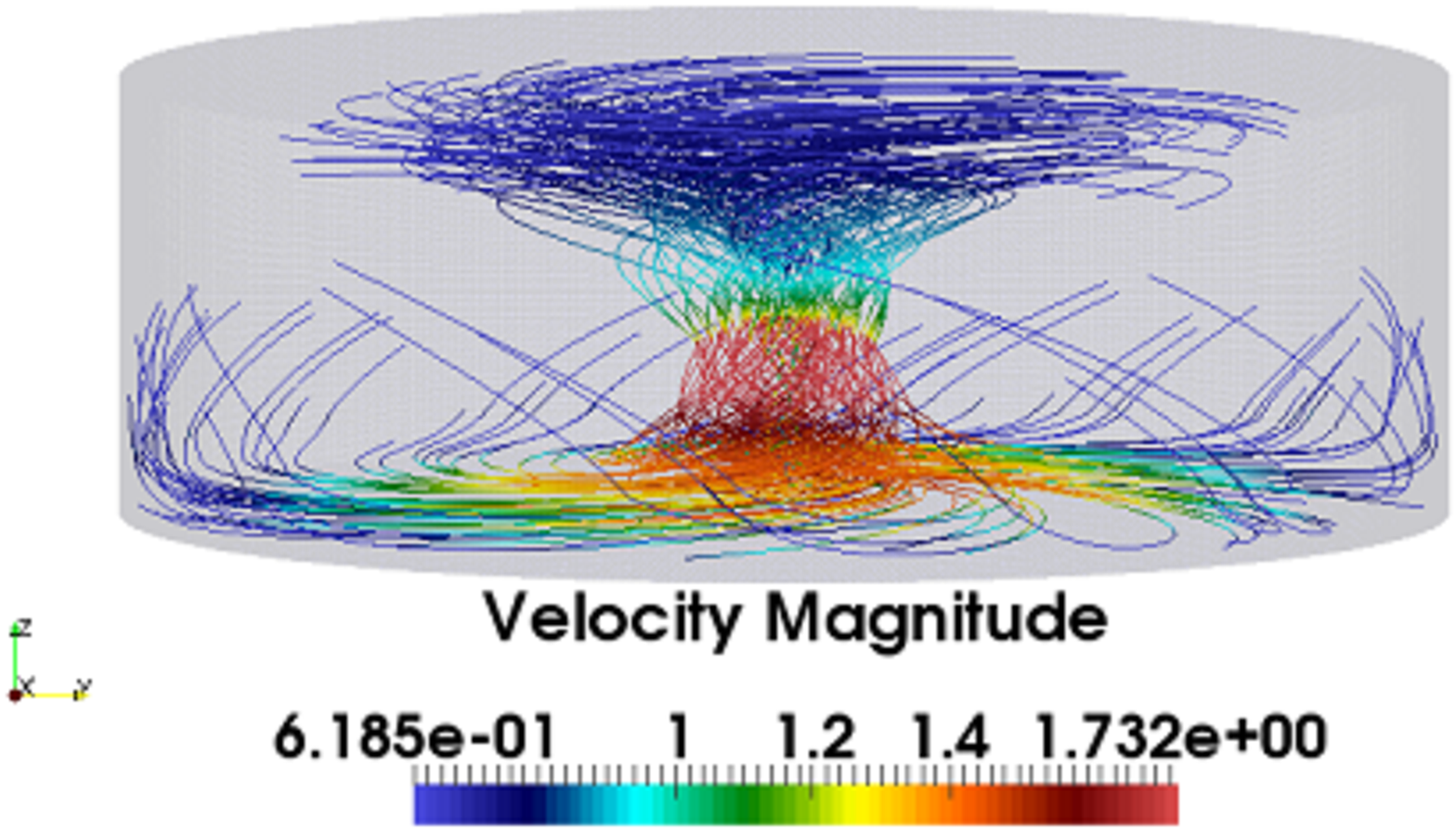}
\includegraphics[width=7.0cm
,keepaspectratio
]{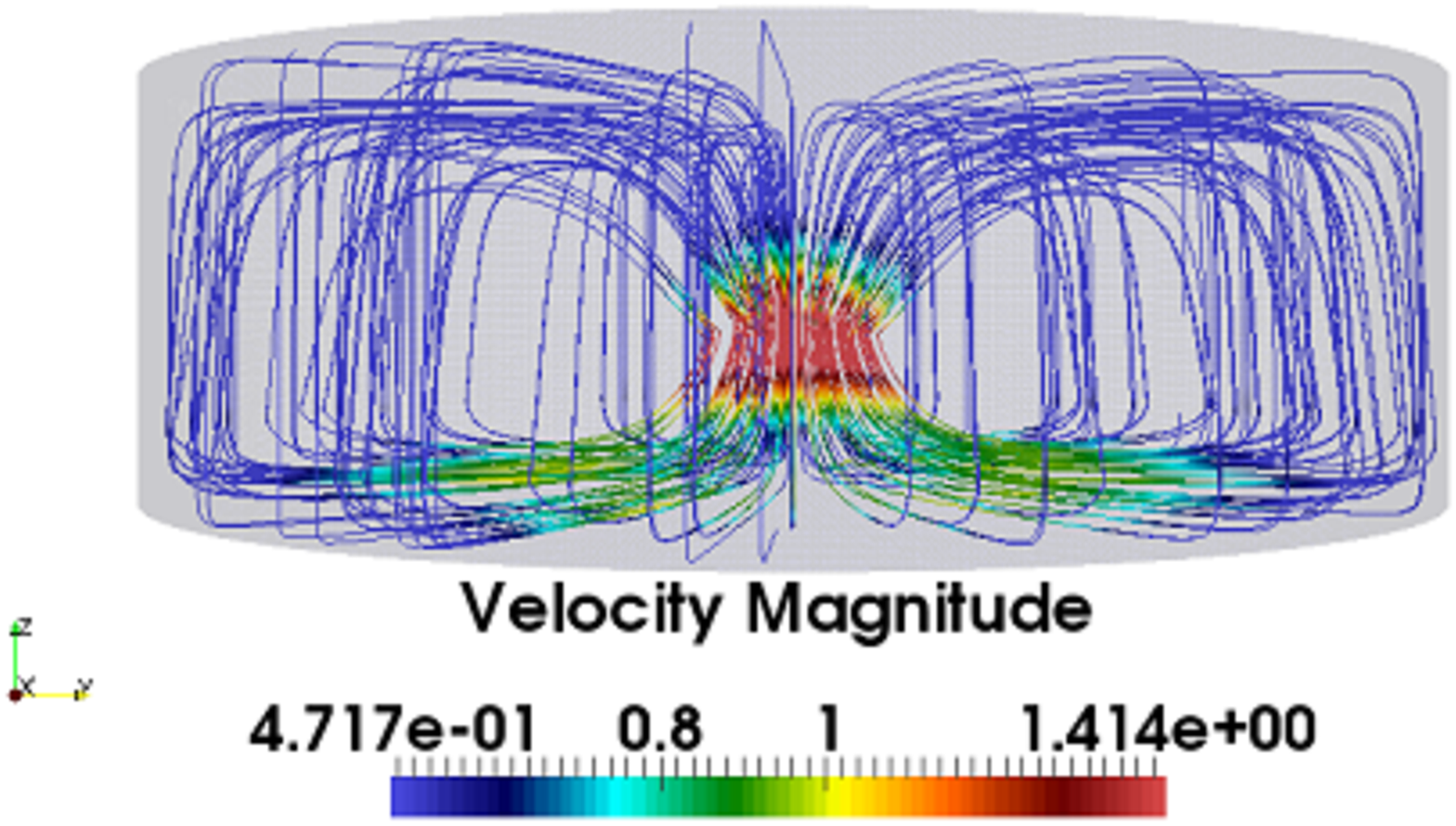}
\vspace{0.5cm}
\caption{Initial data for $a={1}/{8}$, left: swirl, right: no swirl}
\end{figure}


The following figures are our numerical results for $a= 1/8$ and $Re
= 1,000$, $5,000$, $10,000$ and $50,000$.
Hereafter, we call the point attaining the maximum value of
$|v|$ ``the maximum point of $|v|$'' and the cross-section $\{x \in
\Omega; x_1 = 0\}$ ``the plane $x_1=0$'' simply.

\begin{figure}\label{swirlmax}
\includegraphics[width=7.0cm
,keepaspectratio
]{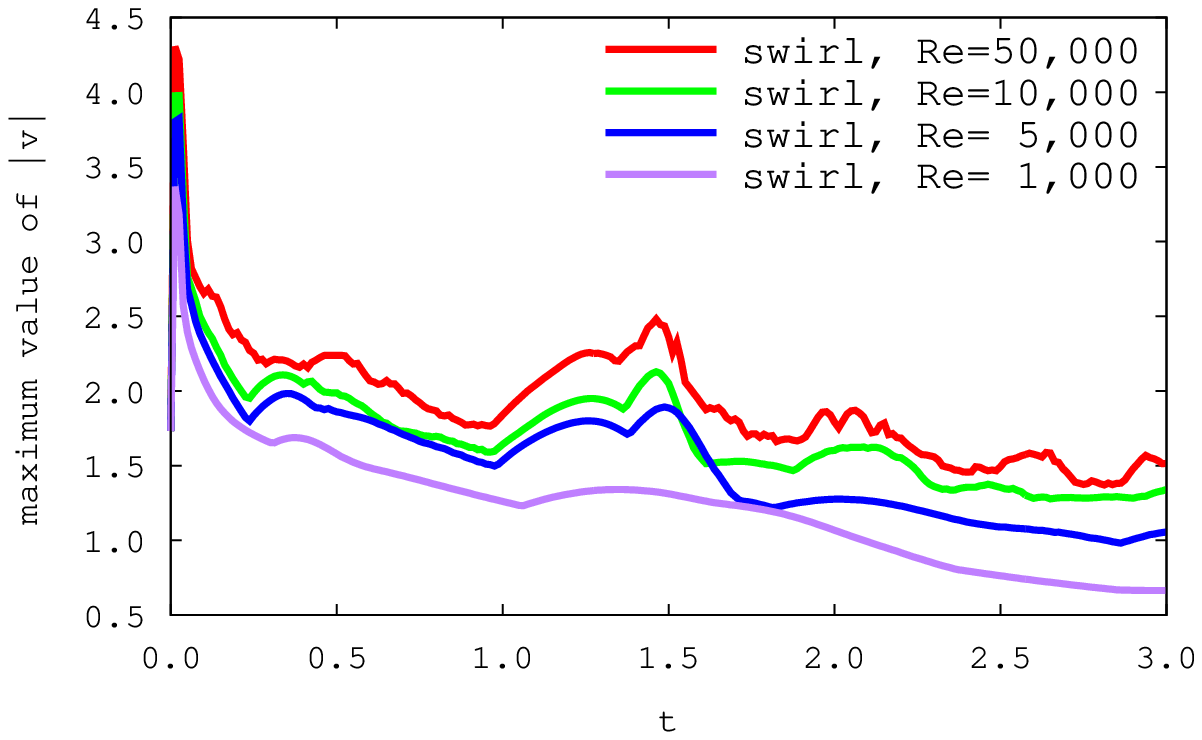}
\includegraphics[width=7.0cm
,keepaspectratio
]{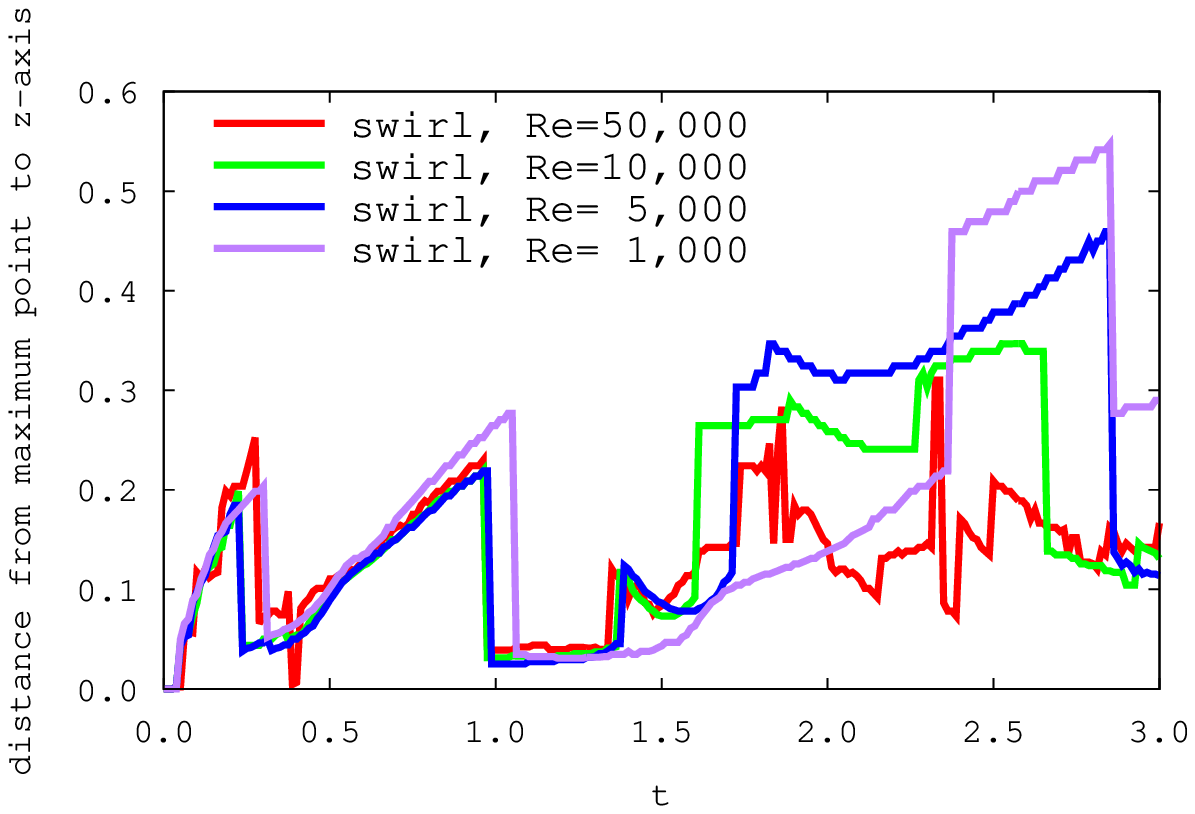}
\vspace{0.5cm}
\caption{Graphs of maximum values of $|v|$ versus $t$ (left) and the
distance from the maximum point of $|v|$ to the $z$-axis versus $t$
(right) for the swirl case with Reynolds numbers $50,000$ (red),
$10,000$ (green), $5,000$ (blue) and $1,000$ (purple).}
\end{figure}
\begin{figure}\label{noswirlmax}
\includegraphics[width=7.0cm
,keepaspectratio
]{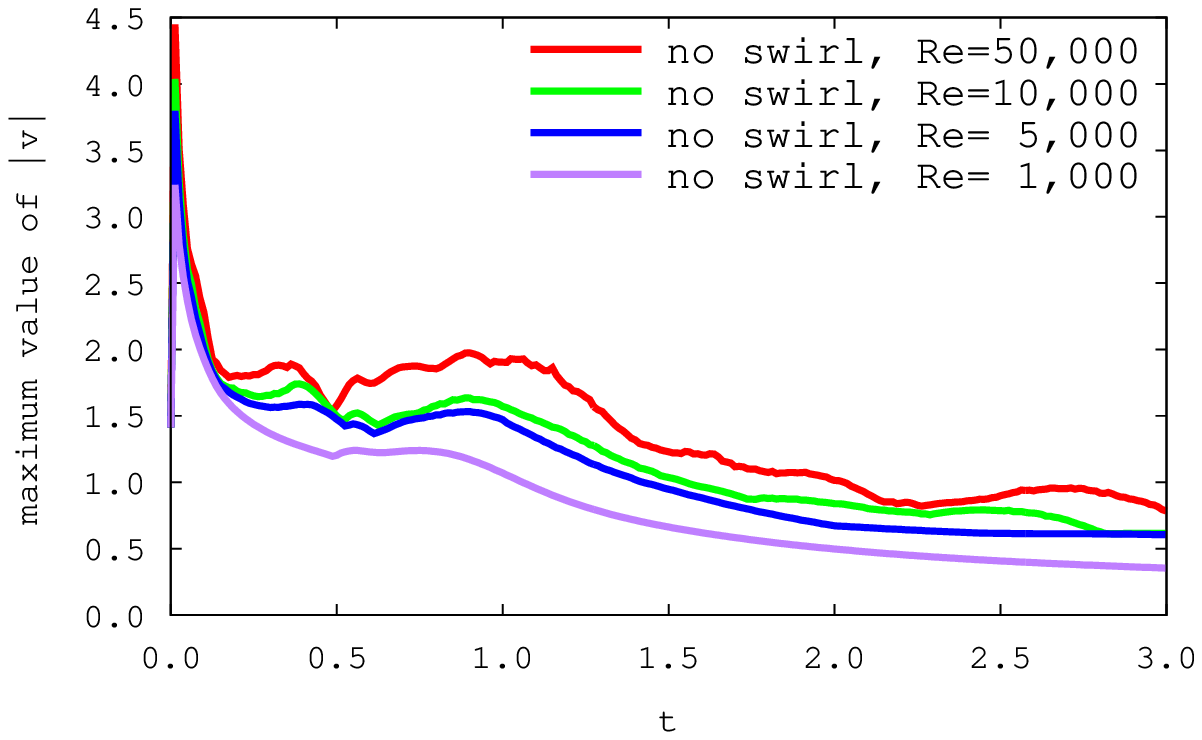}
\includegraphics[width=7.0cm
,keepaspectratio
]{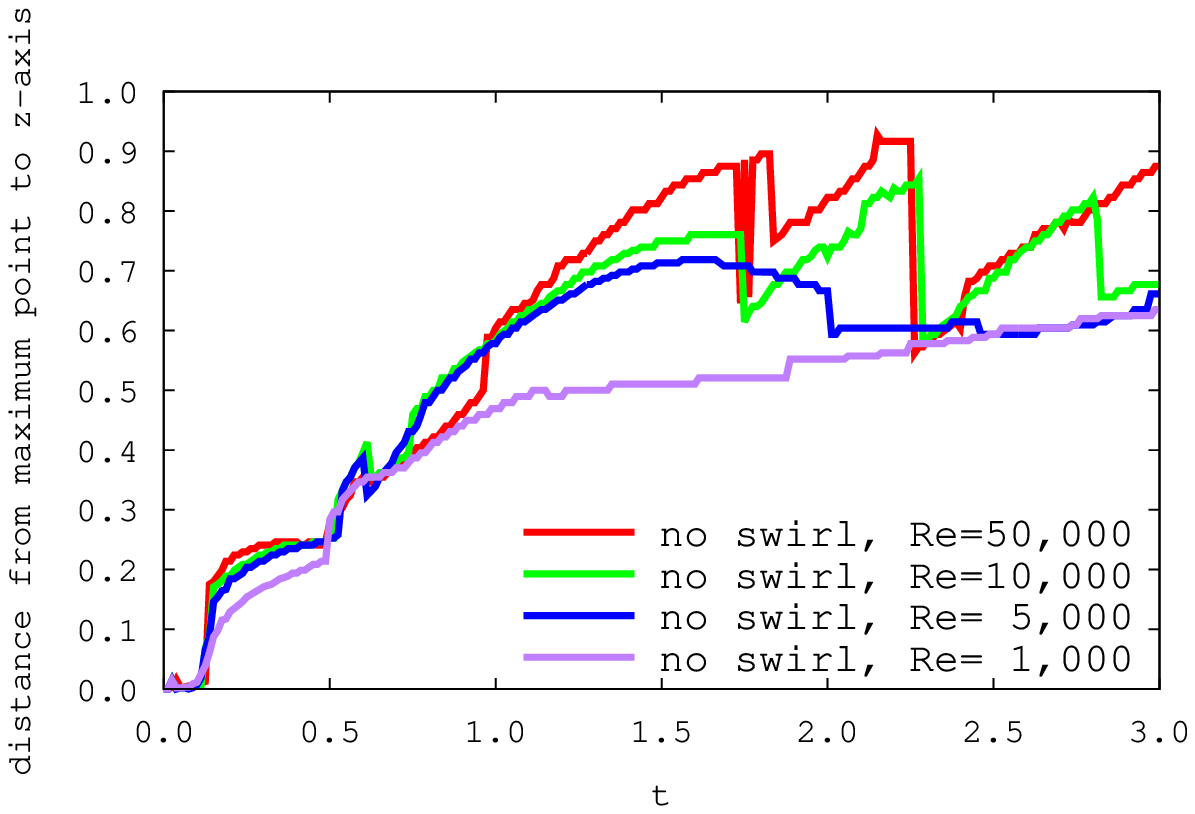}
\vspace{0.5cm}
\caption{Graphs of maximum values of $|v|$ versus $t$ (left) and the
distance from the maximum point of $|v|$ to the $z$-axis versus $t$
(right) in the no swirl case with Reynolds numbers $50,000$ (red),
$10,000$ (green), $5,000$ (blue) and $1,000$ (purple).}
\end{figure}

Figure 2 shows graphs of maximum values of $|v|$ versus time $t$
(left) and the distance from the maximum point of $|v|$ to the
$z$-axis versus time $t$ (right) for the swirl case, where four colors
are used for the graphs of $Re = 50,000$ (red), $10,000$ (green),
$5,000$ (blue) and $1,000$ (purple).
Figure 3 expresses corresponding graphs for the no swirl case.
From figures 2 and 3 we can see that different phenomena appear in the
swirl and no swirl cases.
Typical phenomenon is the drastic changes of the distance from the
maximum point to the z-axis around $t = 0.35$ and $1.0$ in the swirl
case, which does not appear in the no swirl case. 
The difference becomes clearer in a higher Reynolds number flow. In order to understand the difference we display time evolutions
of velocity on the plane $x_1 = 0$ for both cases with $Re = 50,000$
in figures 4, 5 and 6.
Figures 4 and 5 show the time evolutions of $|v|$ on the plane for the
no swirl and swirl cases, respectively, and figure 6 exhibits the time
evolution of the axial velocity $u_z$ ($z$-component of $v$) in the swirl case.
They imply that the flow dissipates straightforwardly as $t$ increases
in the no swirl case and that an interesting flow structure appears
near the z-axis in the swirl case.


\begin{figure}\label{noswirlslice}
\includegraphics[width=7.0cm
,keepaspectratio
]{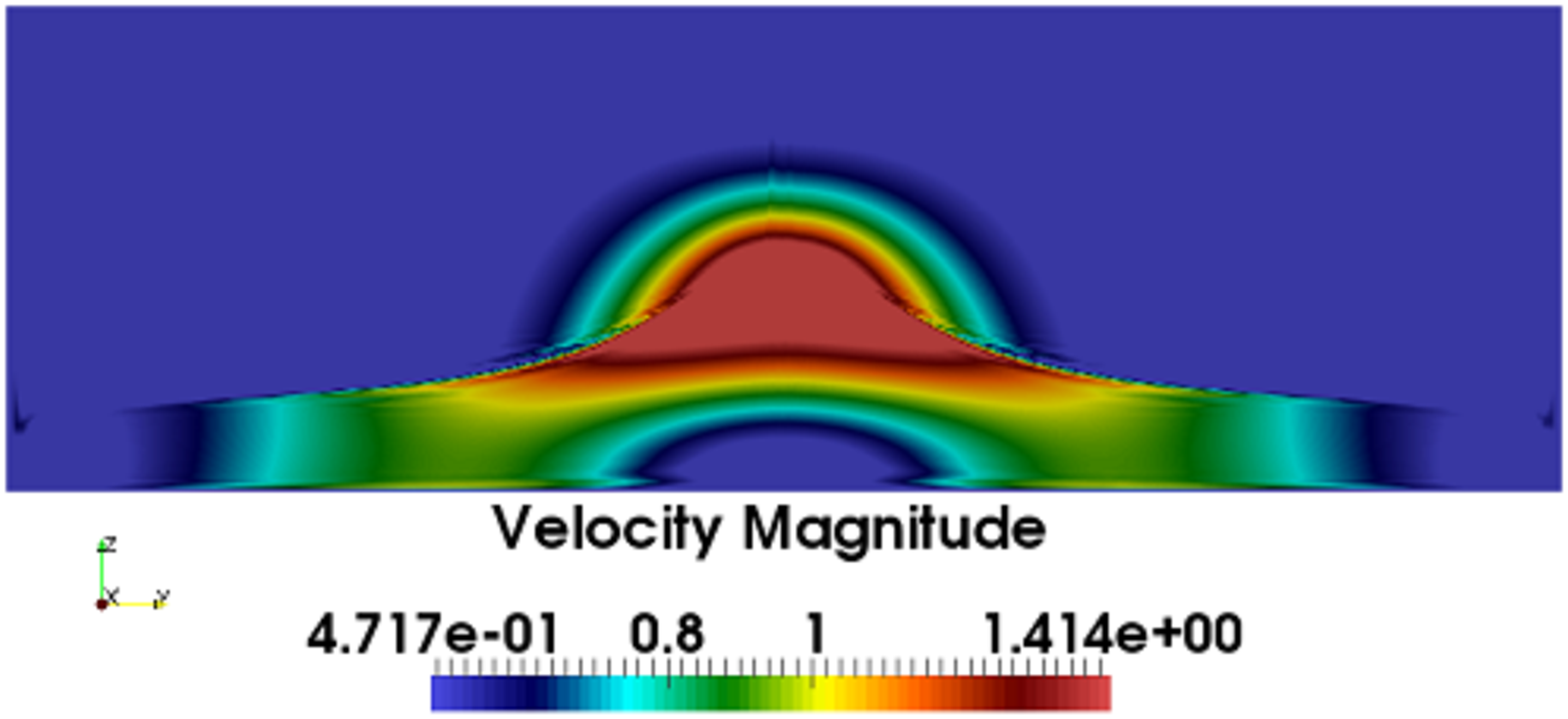}
\includegraphics[width=7.0cm
,keepaspectratio
]{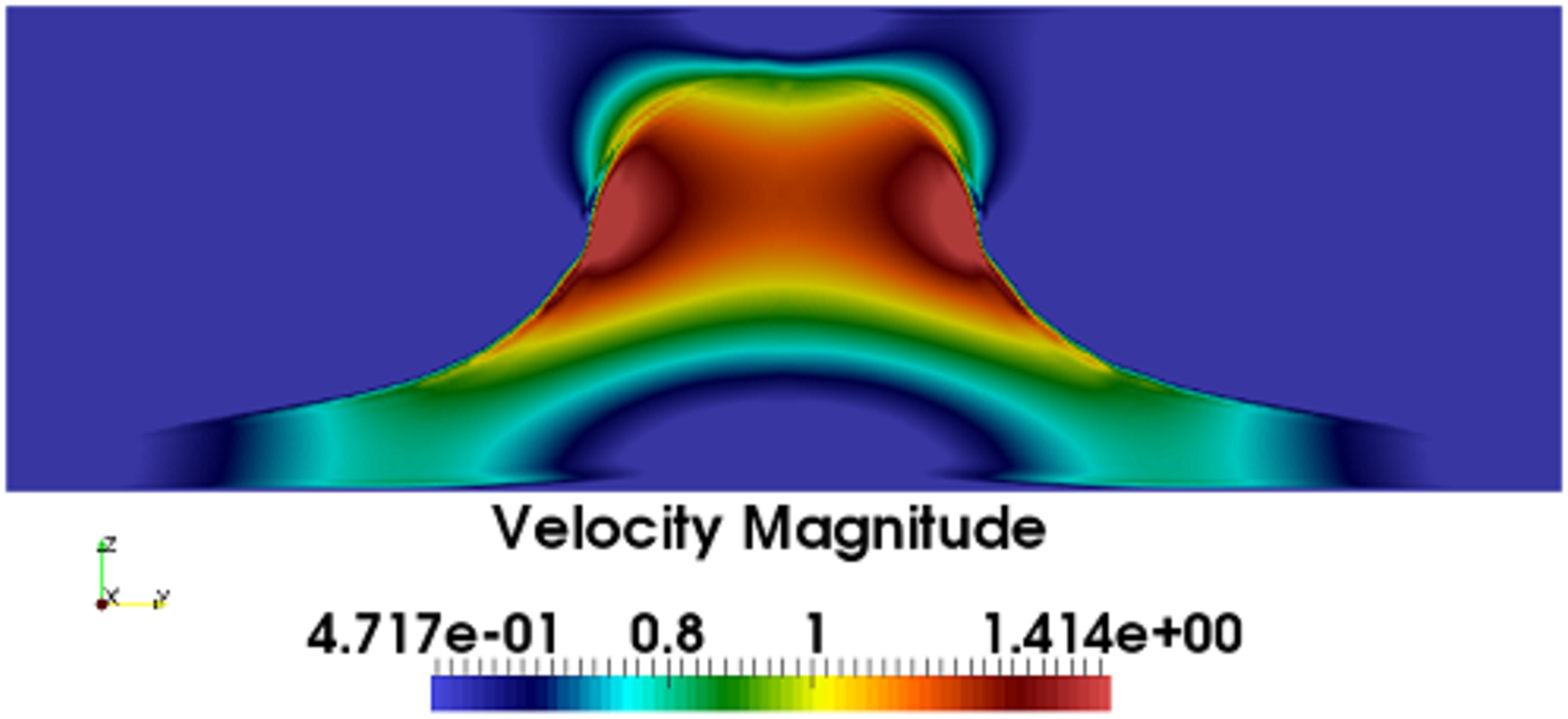}
\includegraphics[width=7.0cm
,keepaspectratio
]{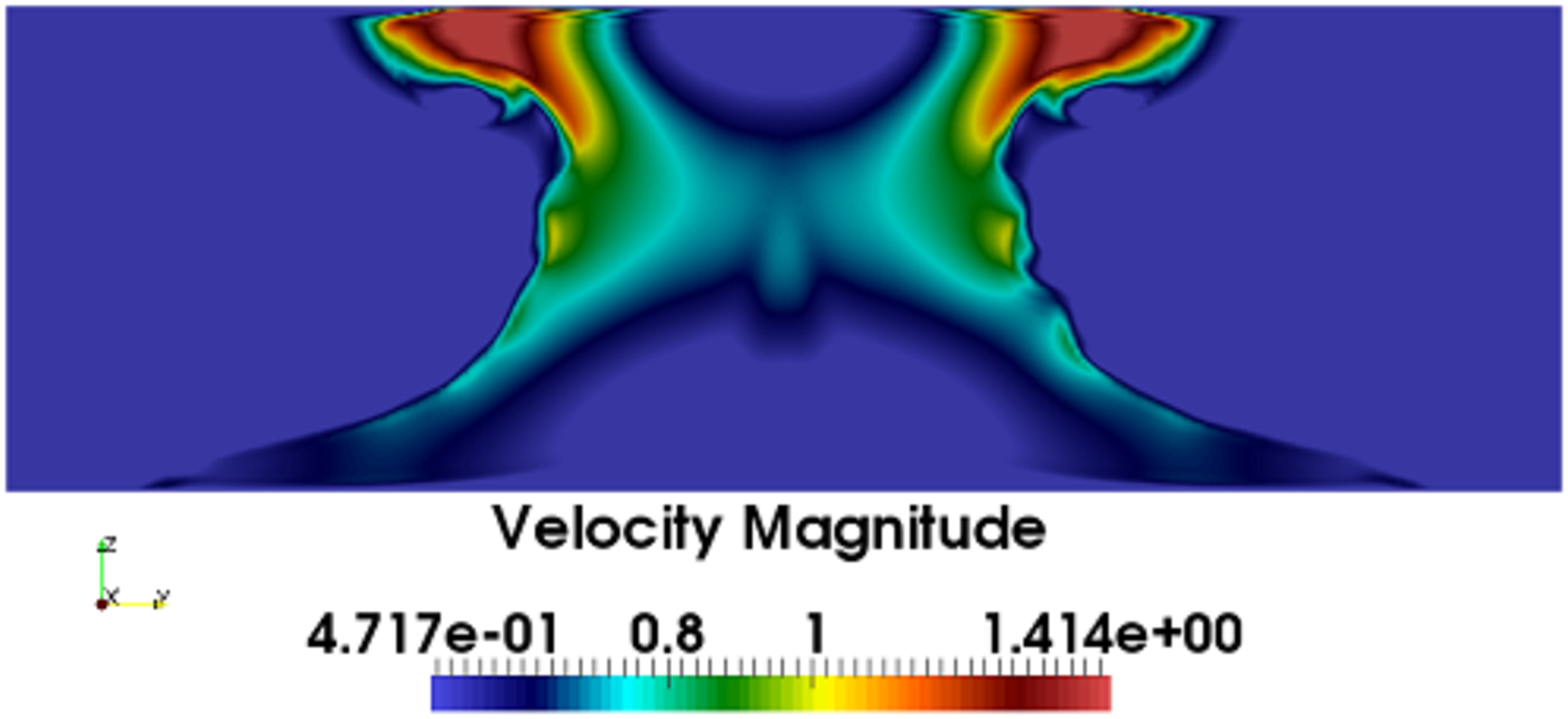}
\includegraphics[width=7.0cm
,keepaspectratio
]{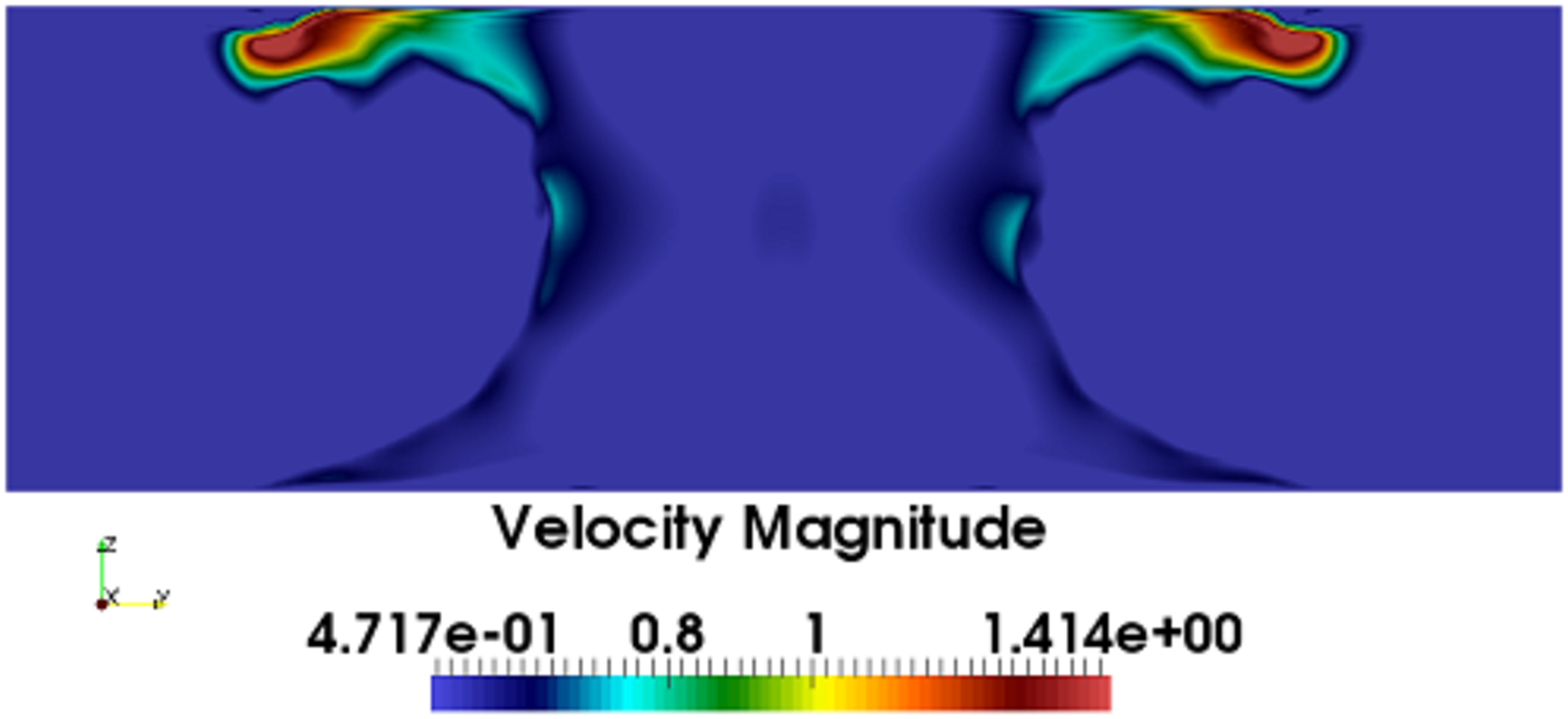}
\includegraphics[width=7.0cm
,keepaspectratio
]{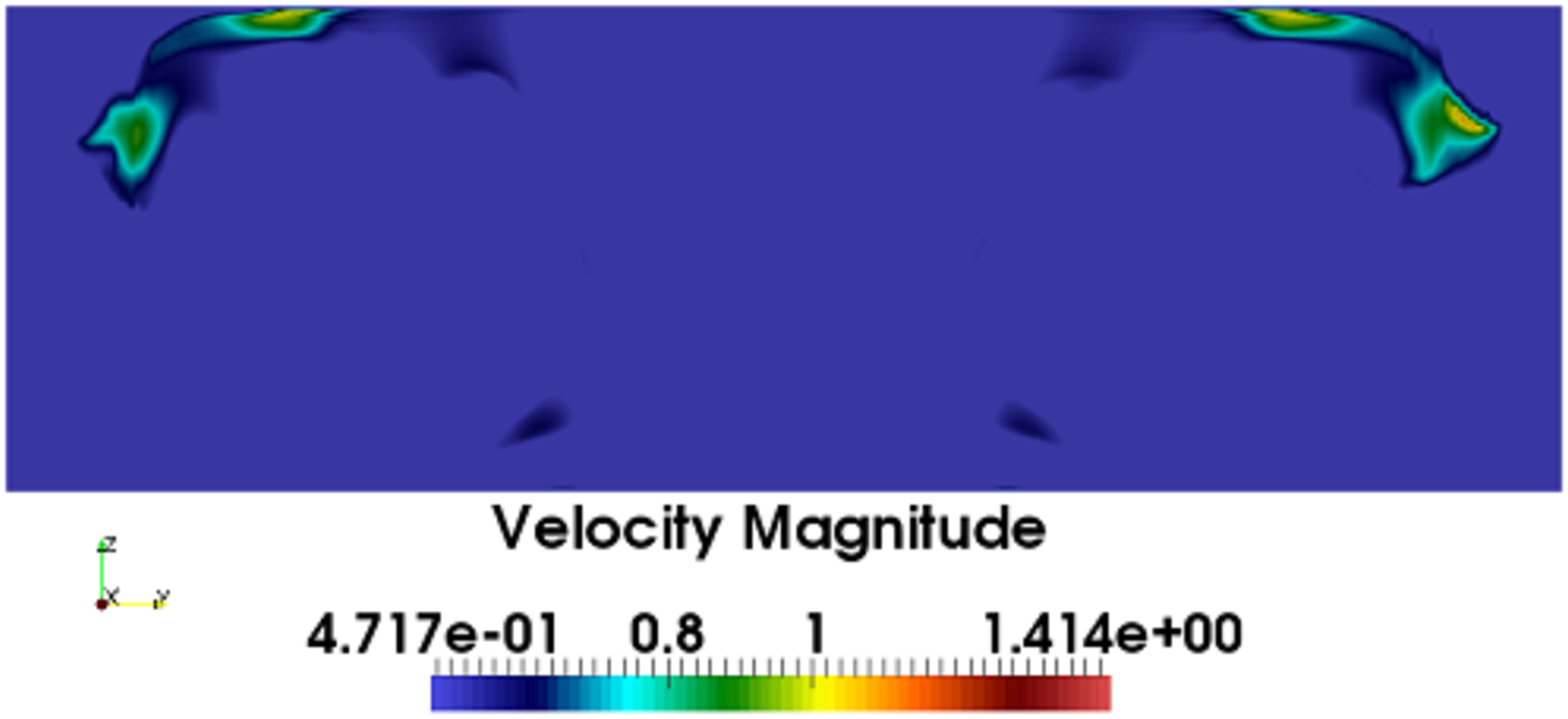}
\includegraphics[width=7.0cm
,keepaspectratio
]{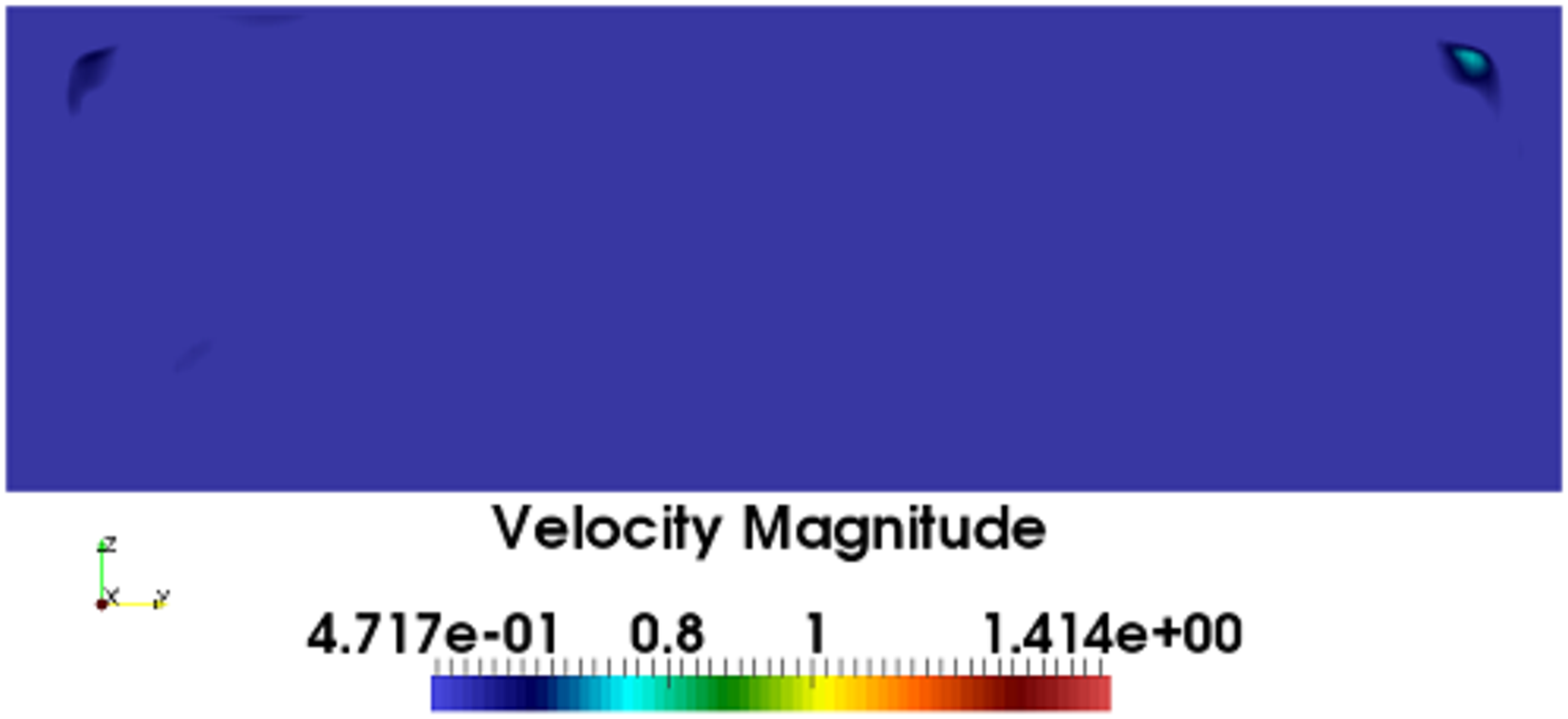}
\vspace{0.5cm}
\caption{Time evolution of the velocity magnitude $|v|$ on the plane
$x_1=0$ in the no swirl case with $Re = 50,000$. $t = 0.1$ (top left),
$0.3$ (top right), $0.7$ (middle left), $1.1$ (middle right), $1.7$
(bottom left) and $3.0$ (bottom right).}
\end{figure}

From figures 2, 5 and 6 it is observed that a downward flow arises
near the $z$-axis around $t = 0.3$, that the maximum value of $|v|$ is
attained near the $z$-axis and the lower boundary around the same time
$(t = 0.3)$, that a new upward flow arises near the z-axis around $t =
1.3$, and that the velocity attains its maximum value near the
$z$-axis and the upper boundary around the same time $(t = 1.3)$. 

\begin{rem}
The scale may be different in each figures. Red region represents high magnitude.
\end{rem}


\begin{figure}\label{swirlslice}
\includegraphics[width=7.0cm
,keepaspectratio
]{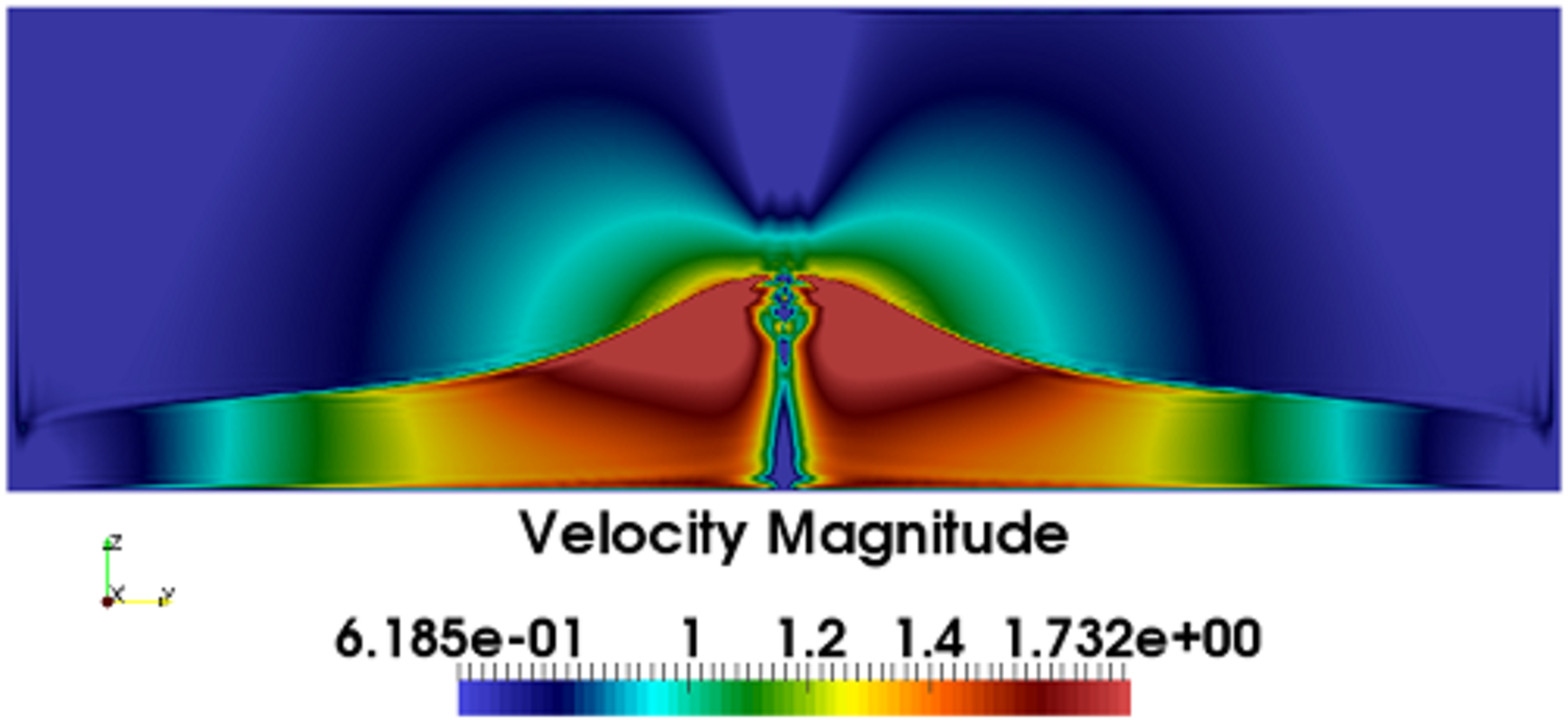}
\includegraphics[width=7.0cm
,keepaspectratio
]{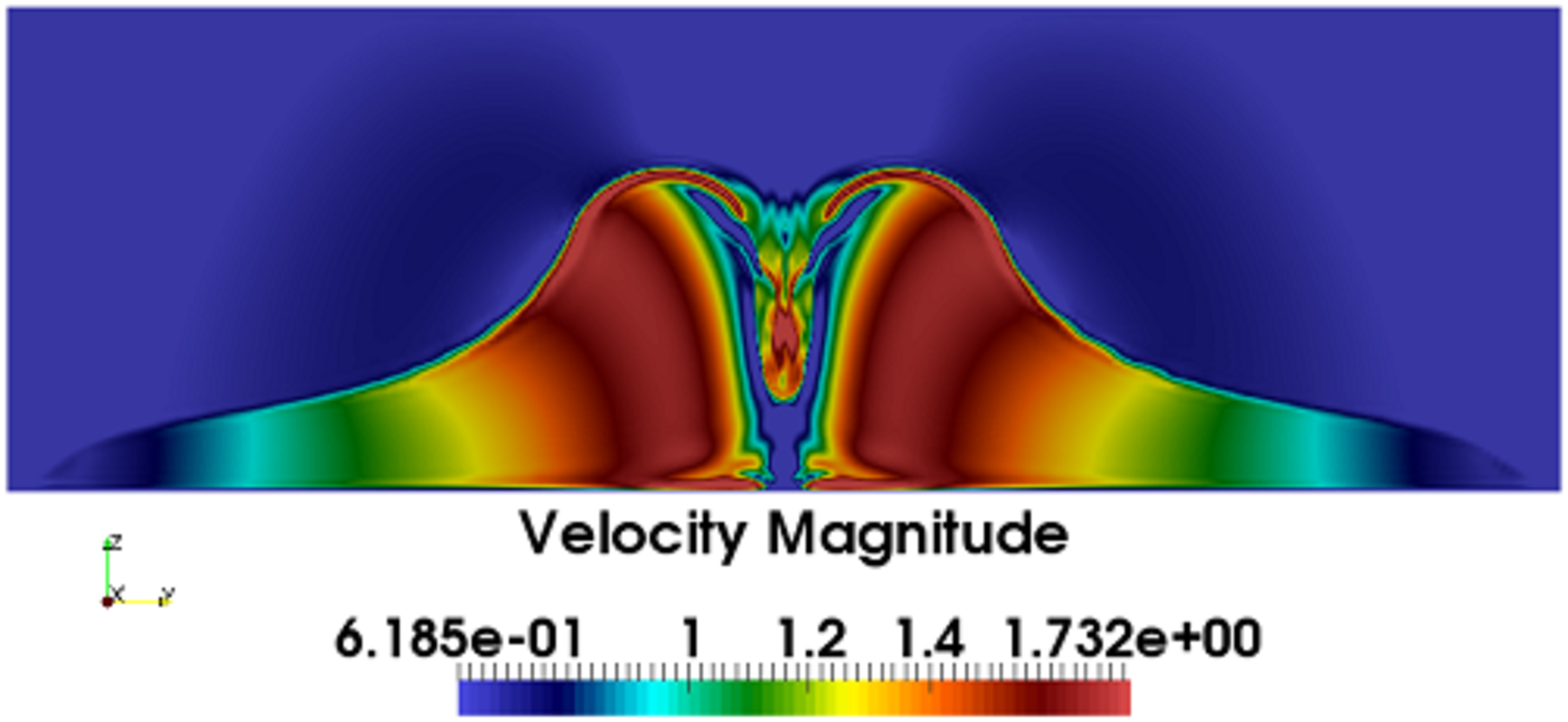}
\includegraphics[width=7.0cm
,keepaspectratio
]{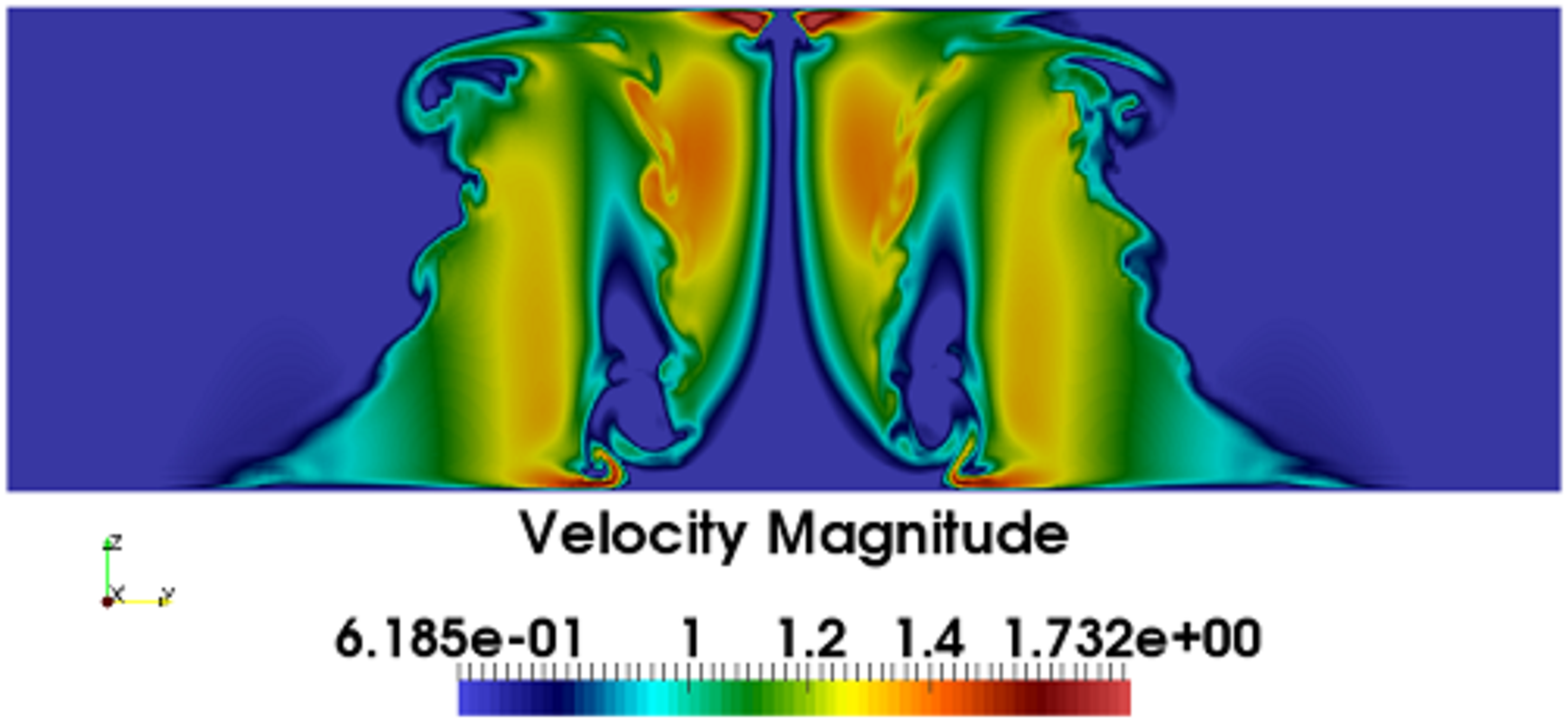}
\includegraphics[width=7.0cm
,keepaspectratio
]{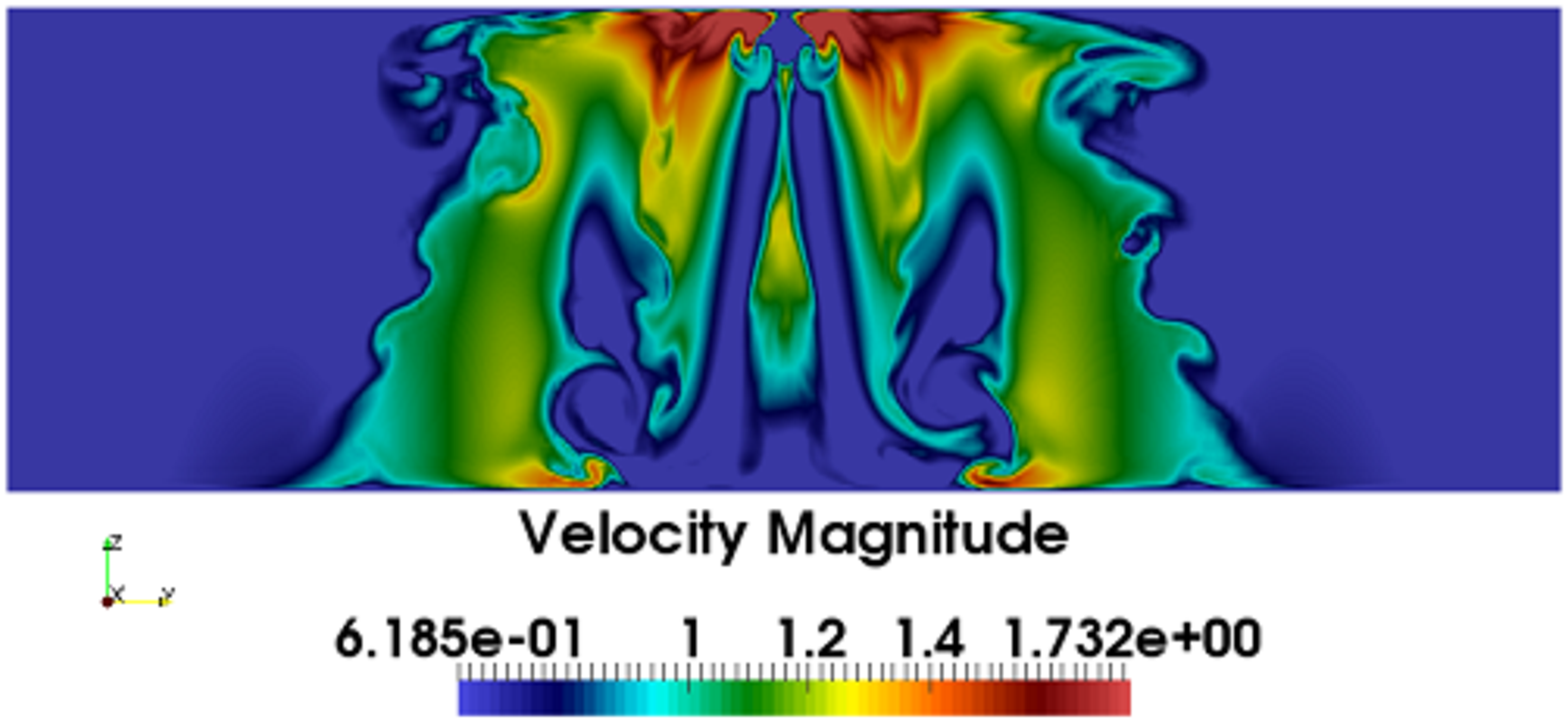}
\includegraphics[width=7.0cm
,keepaspectratio
]{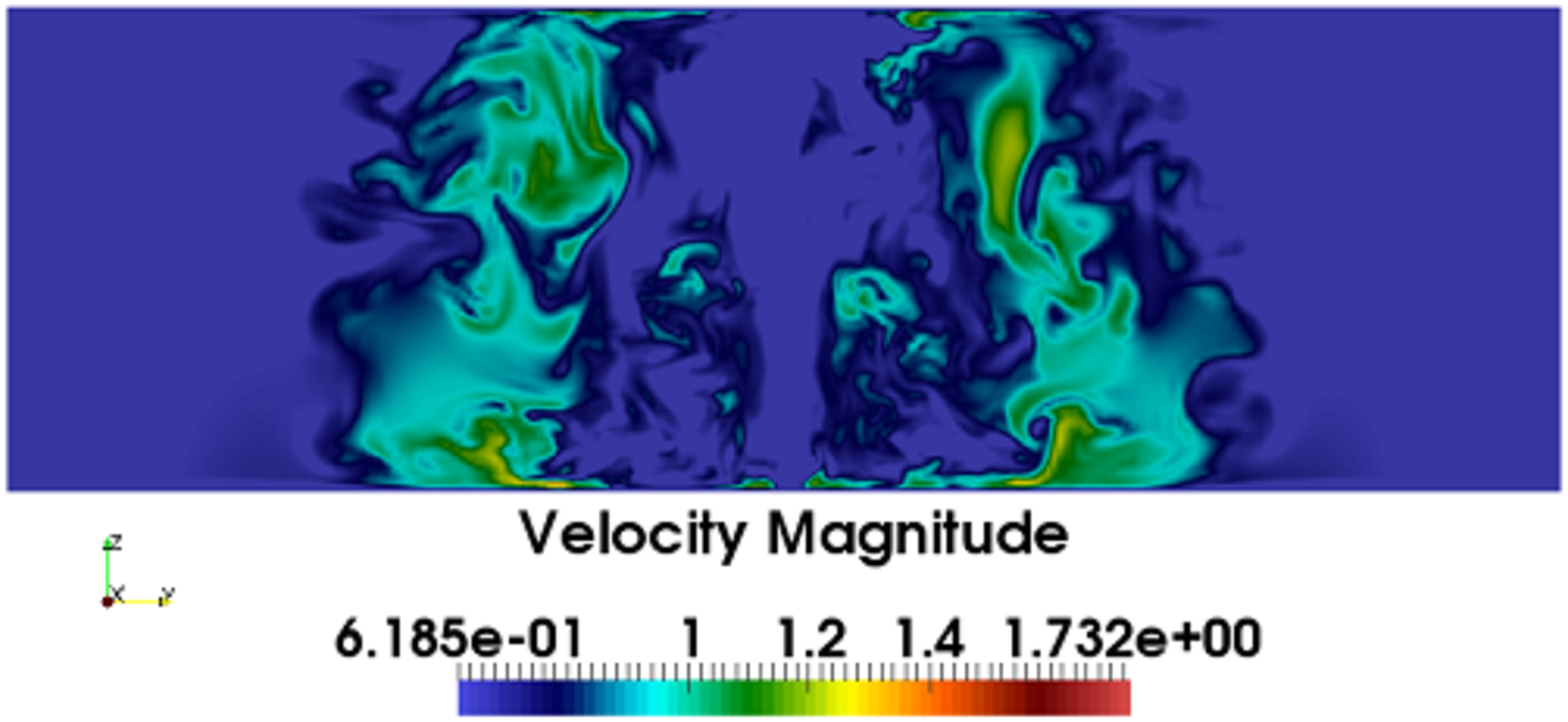}
\includegraphics[width=7.0cm
,keepaspectratio
]{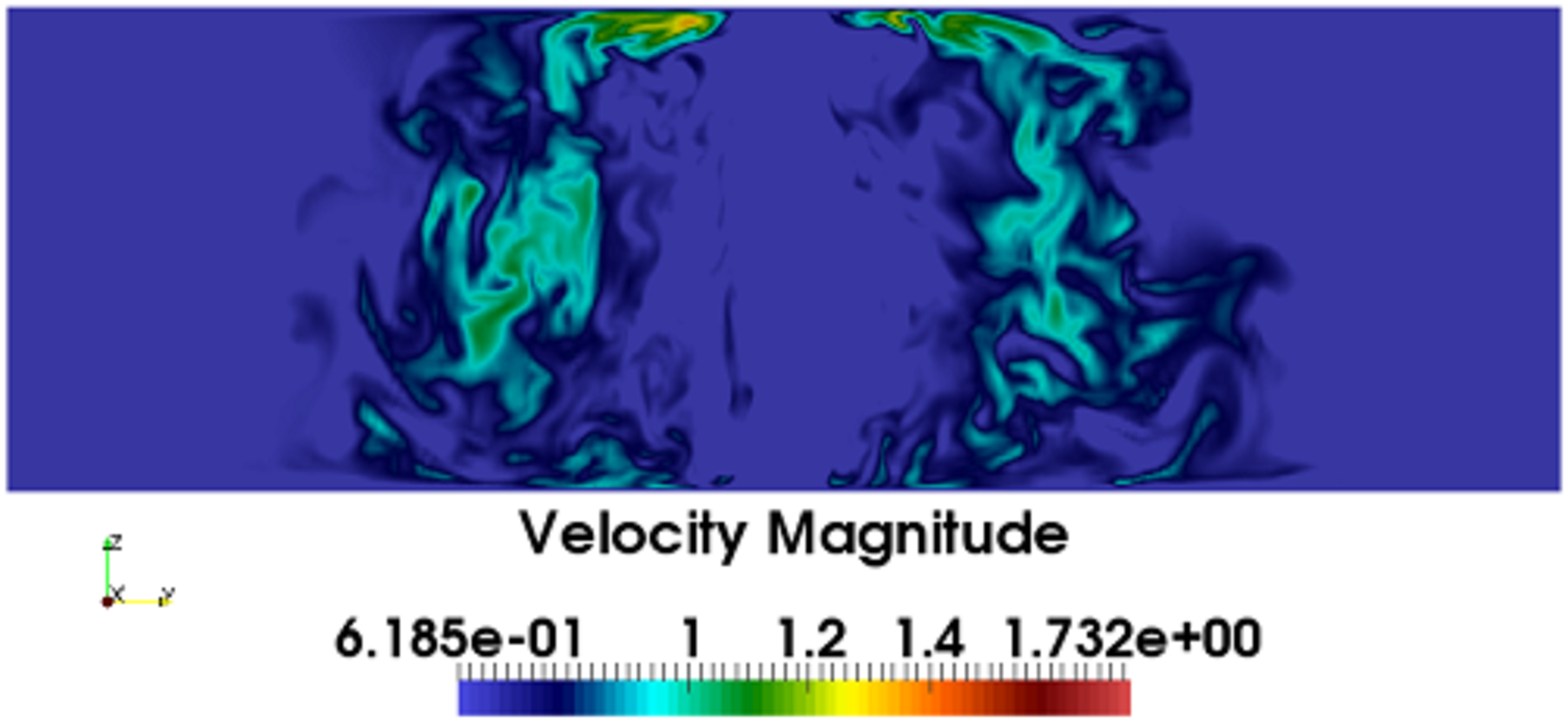}
\vspace{0.5cm}
\caption{Time evolution of the velocity magnitude $|v|$ on the plane
$x_1=0$ in the swirl case with $Re = 50,000$. $t = 0.1$ (top left),
$0.3$ (top right), $1.0$ (middle left), $1.3$ (middle right), $2.3$
(bottom left) and $3.0$ (bottom right).}
\end{figure}
\begin{figure}\label{vzslice}
\includegraphics[width=7.0cm
,keepaspectratio
]{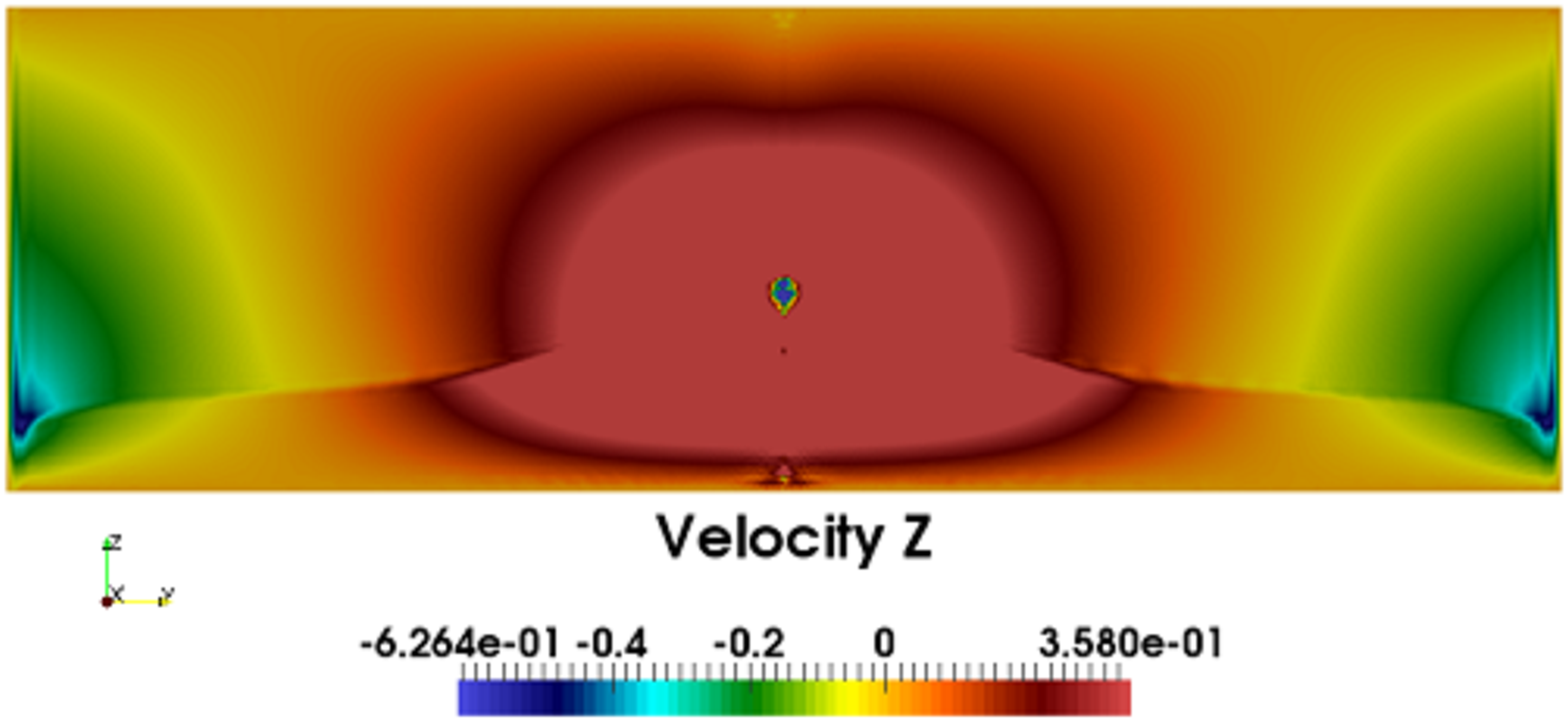}
\includegraphics[width=7.0cm
,keepaspectratio
]{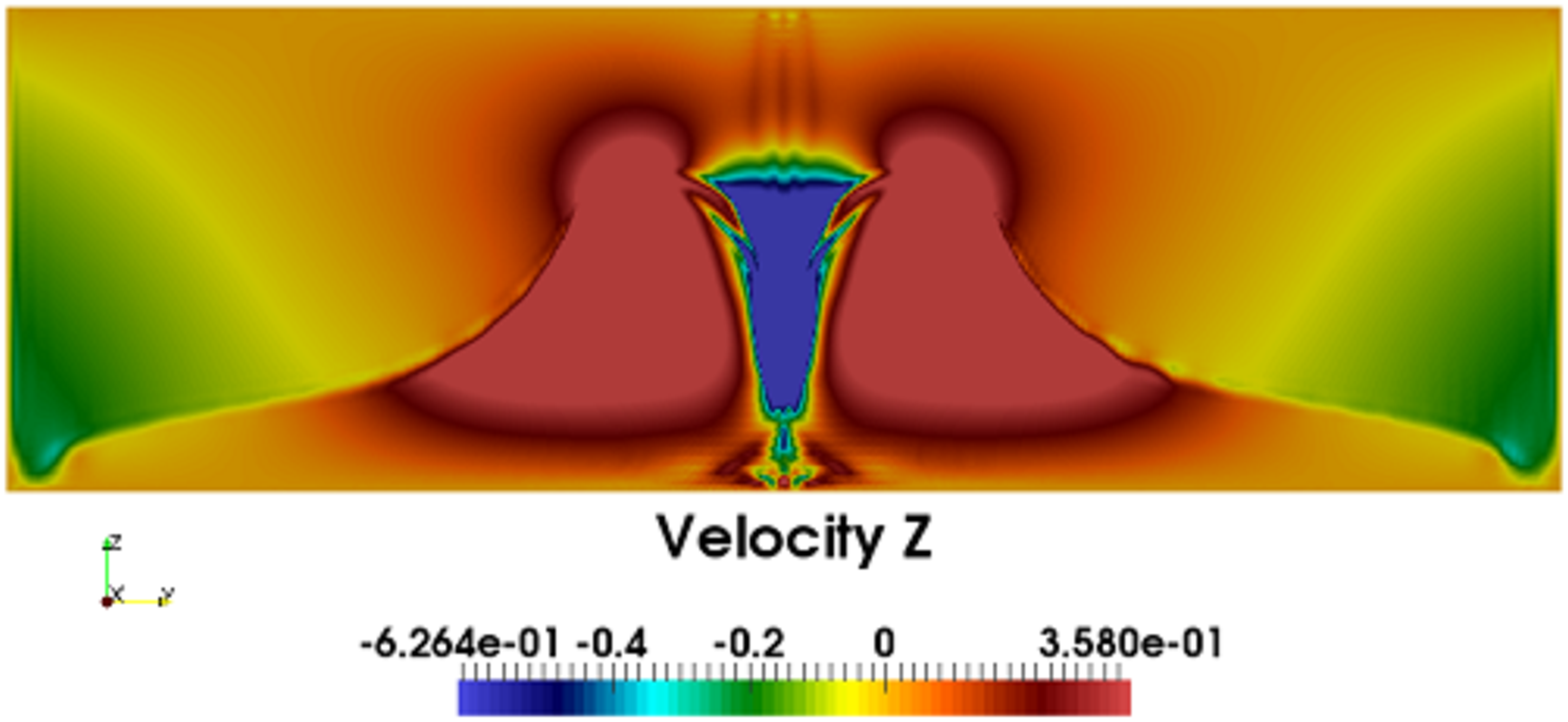}
\includegraphics[width=7.0cm
,keepaspectratio
]{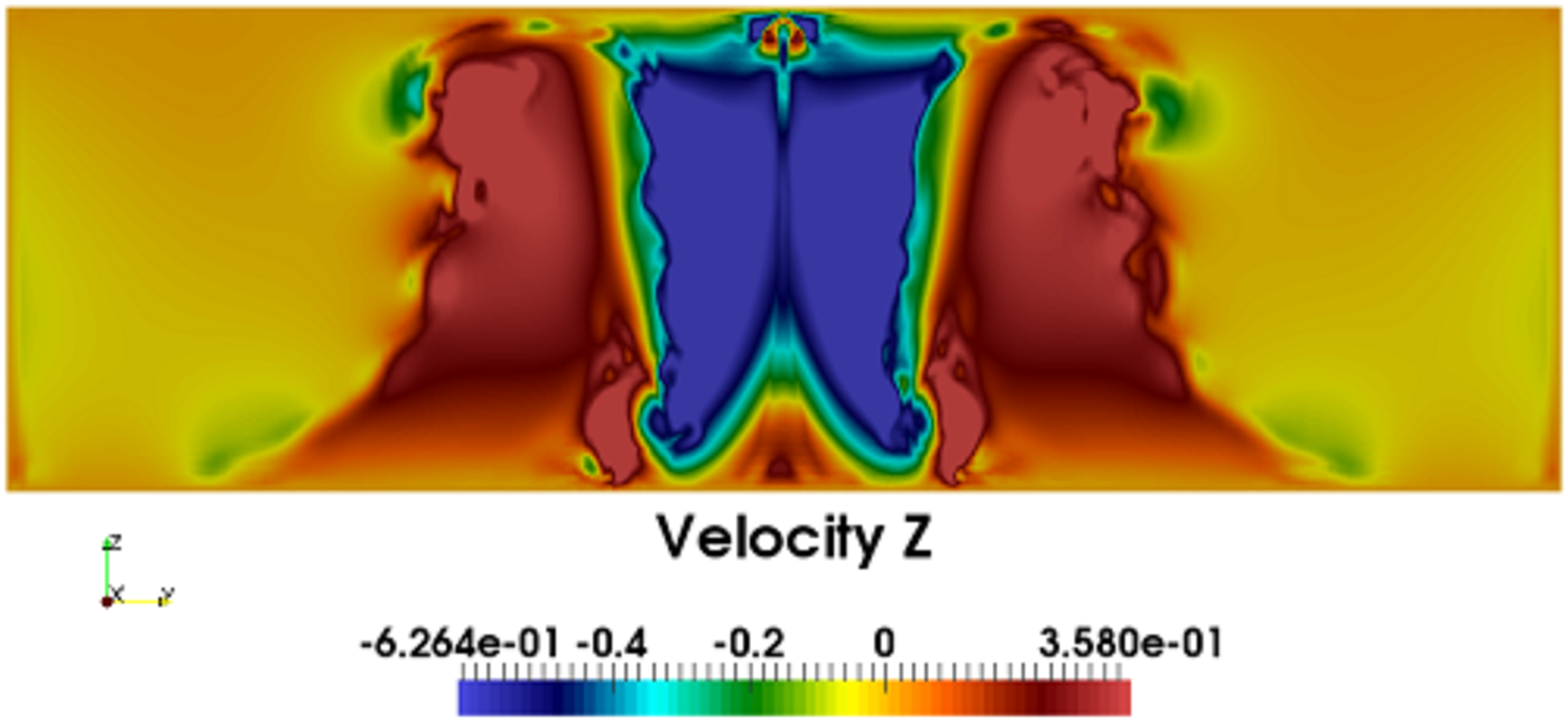}
\includegraphics[width=7.0cm
,keepaspectratio
]{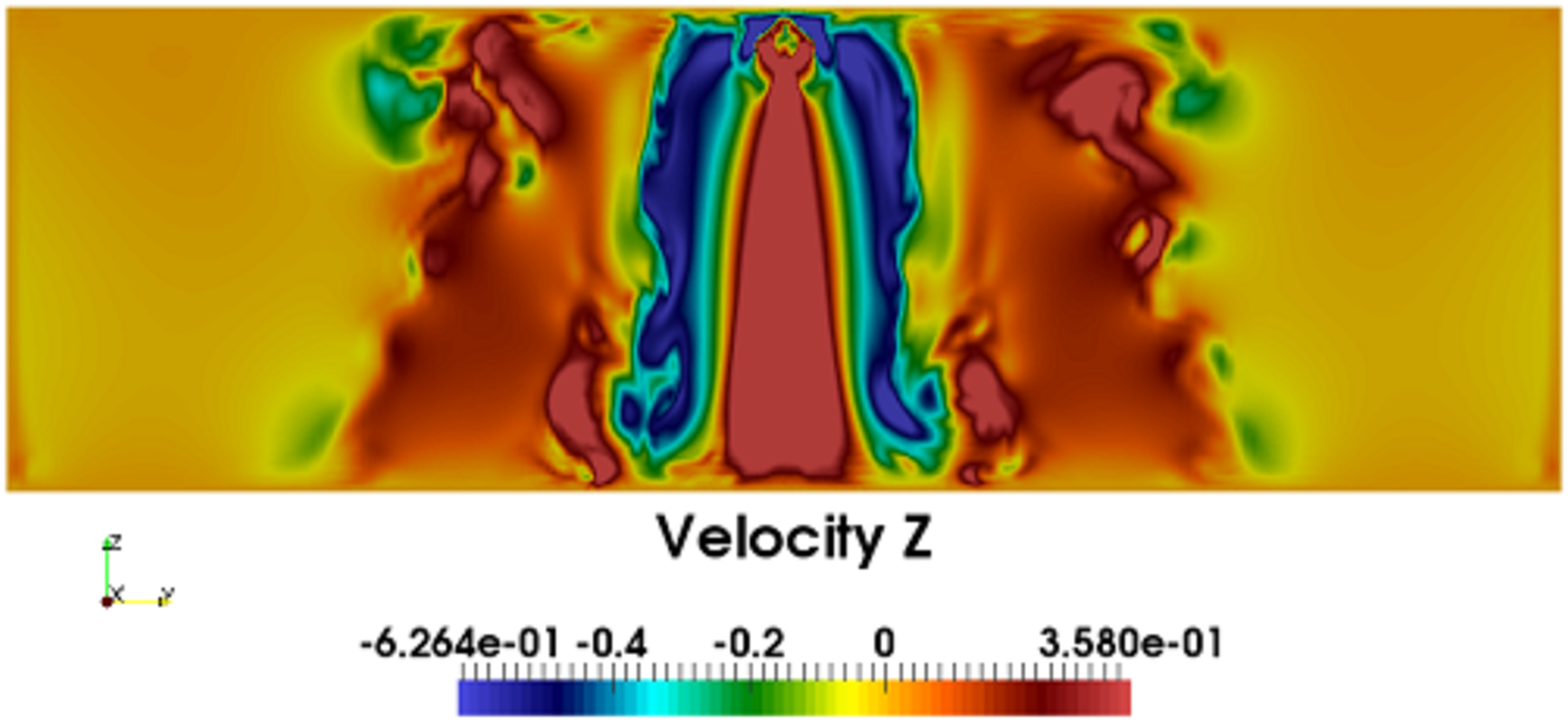}
\includegraphics[width=7.0cm
,keepaspectratio
]{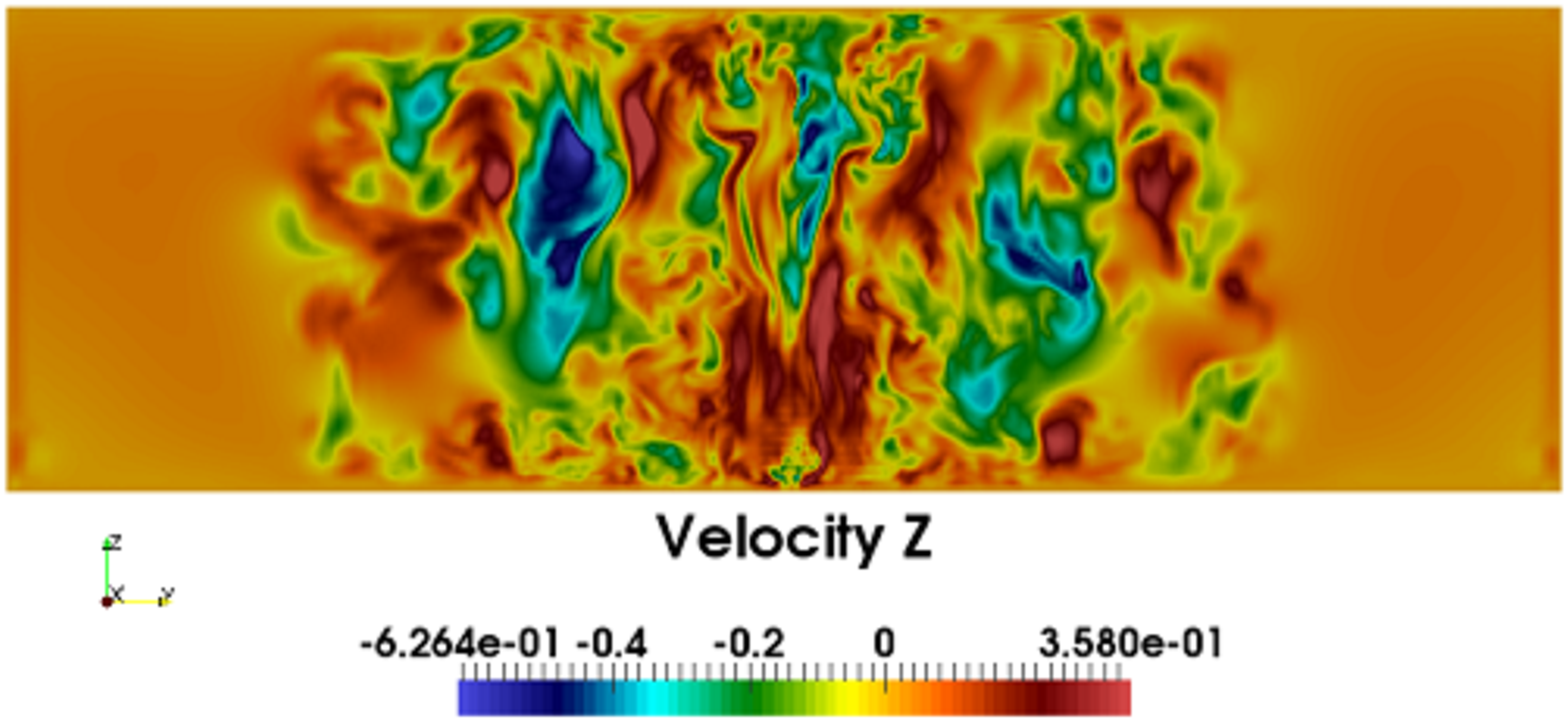}
\includegraphics[width=7.0cm
,keepaspectratio
]{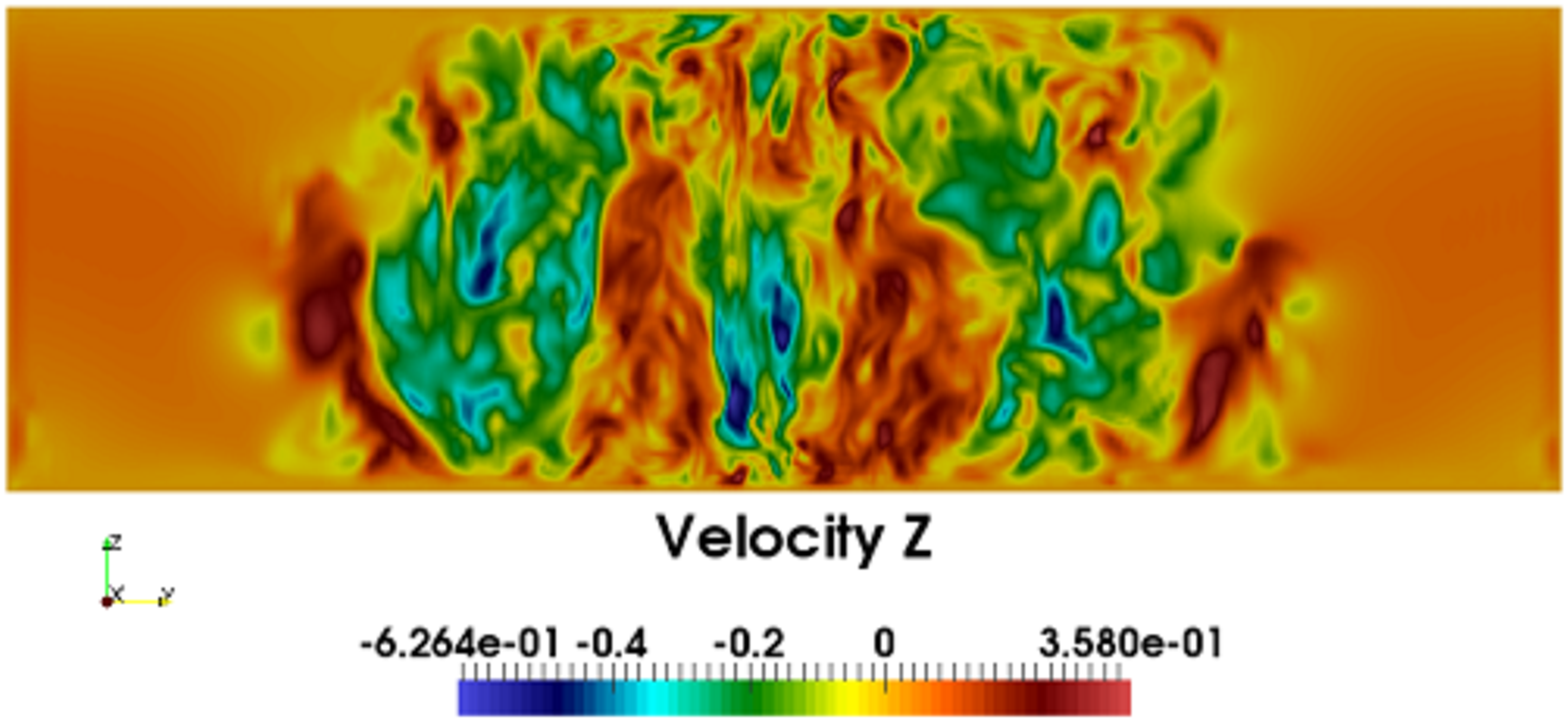}
\vspace{0.5cm}
\caption{Time evolution of the axial velocity $u_z$ on the
plane $x_1=0$ in the swirl case with $Re = 50,000$. $t = 0.1$ (top
left), $0.3$ (top right), $1.0$ (middle left), $1.3$ (middle right),
$2.3$ (bottom left) and $3.0$ (bottom right). Note that red and blue
colors represent the positive and negative values in this figure.}
\end{figure}


\begin{figure}\label{keytime}
\begin{center}
\vspace{-2.1mm}
\includegraphics[width=8.8cm
,keepaspectratio
]{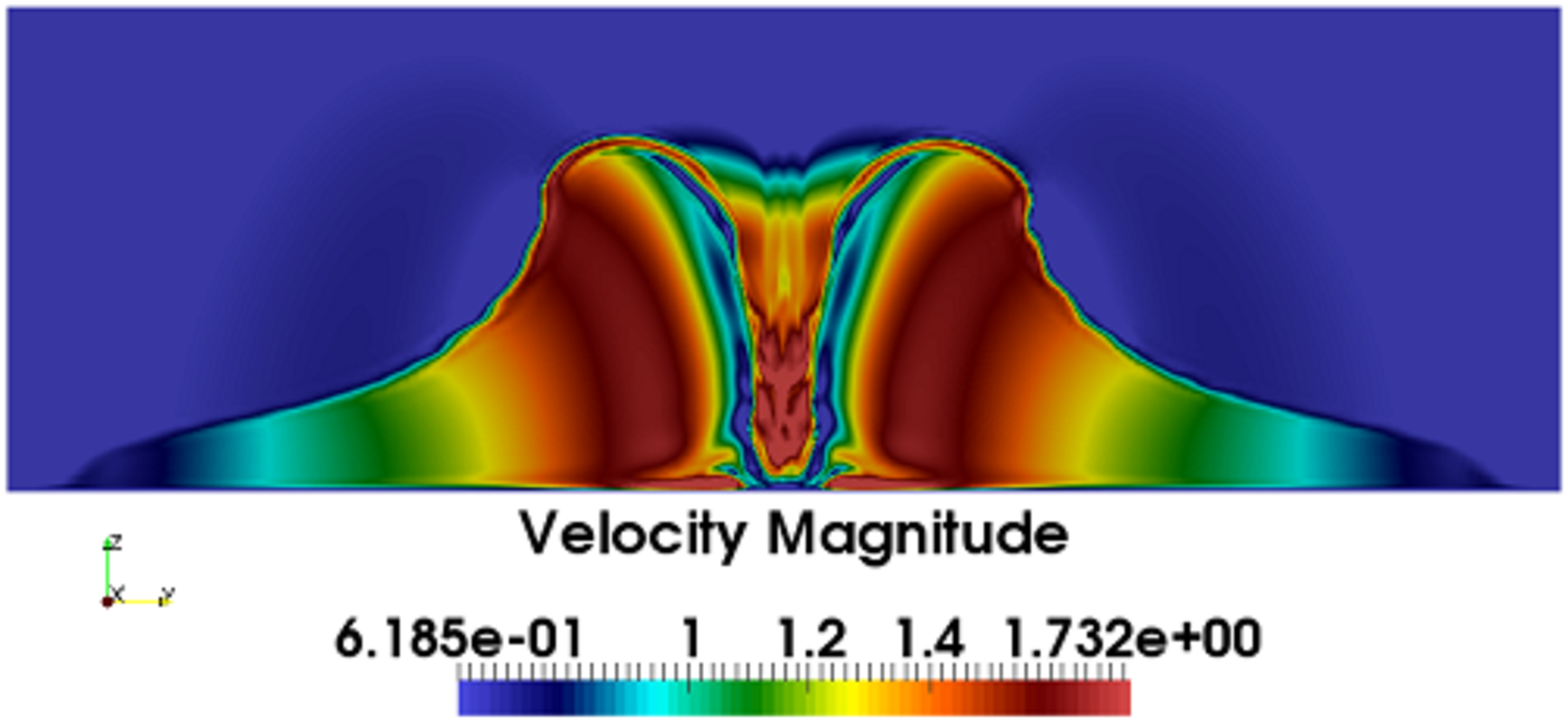}

\vspace{-1.2mm}

\includegraphics[width=8.8cm
,keepaspectratio
]{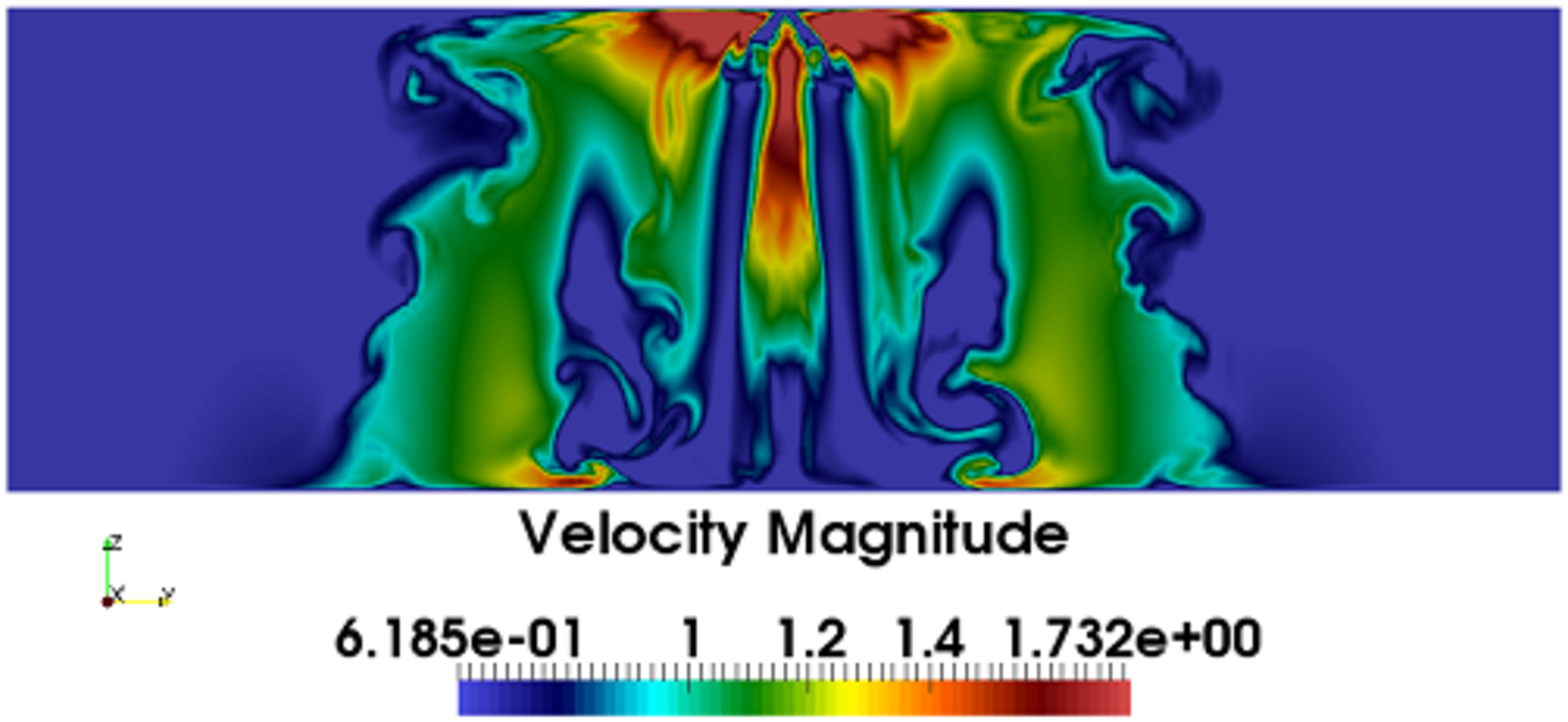}

\vspace{-1.2mm}

\includegraphics[width=8.8cm
,keepaspectratio
]{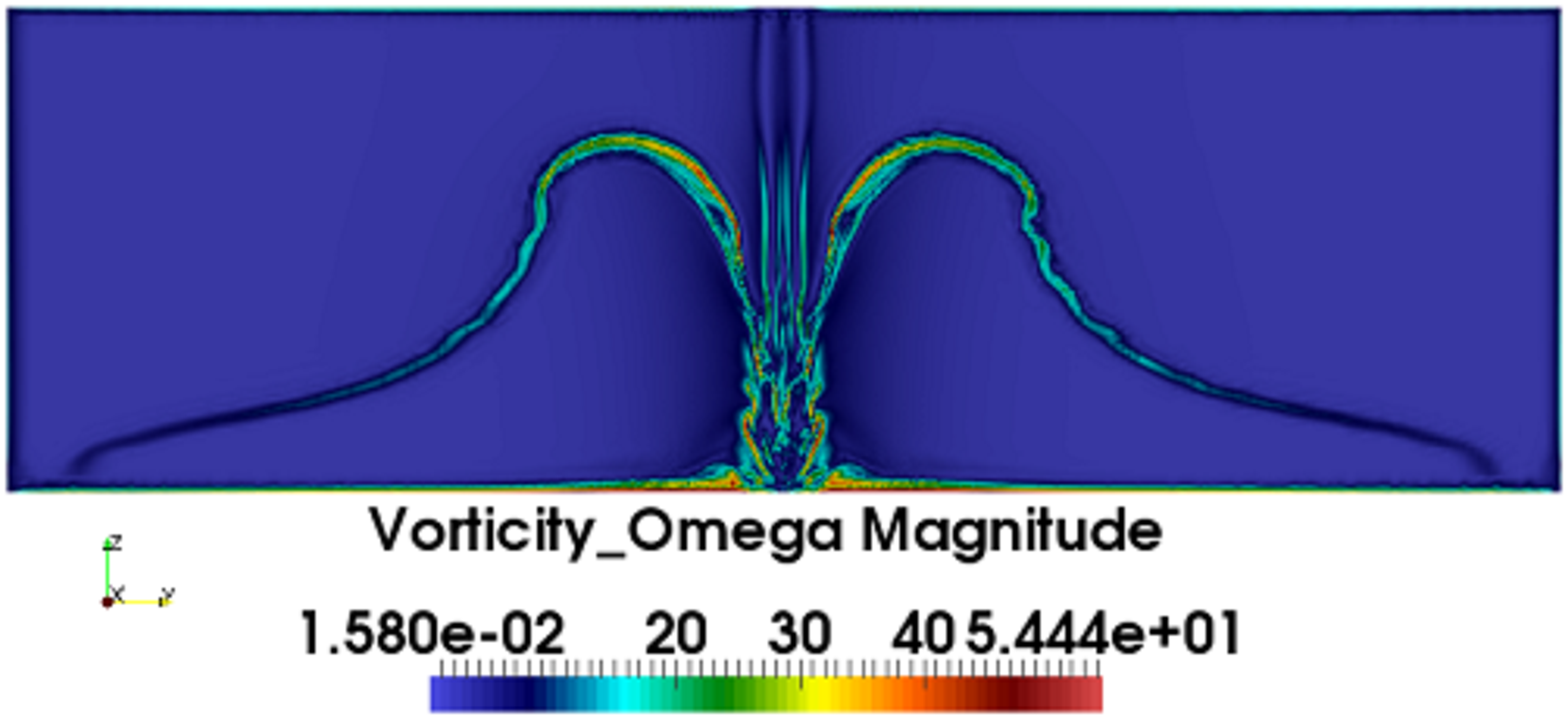}

\vspace{-1.2mm}

\includegraphics[width=8.8cm
,keepaspectratio
]{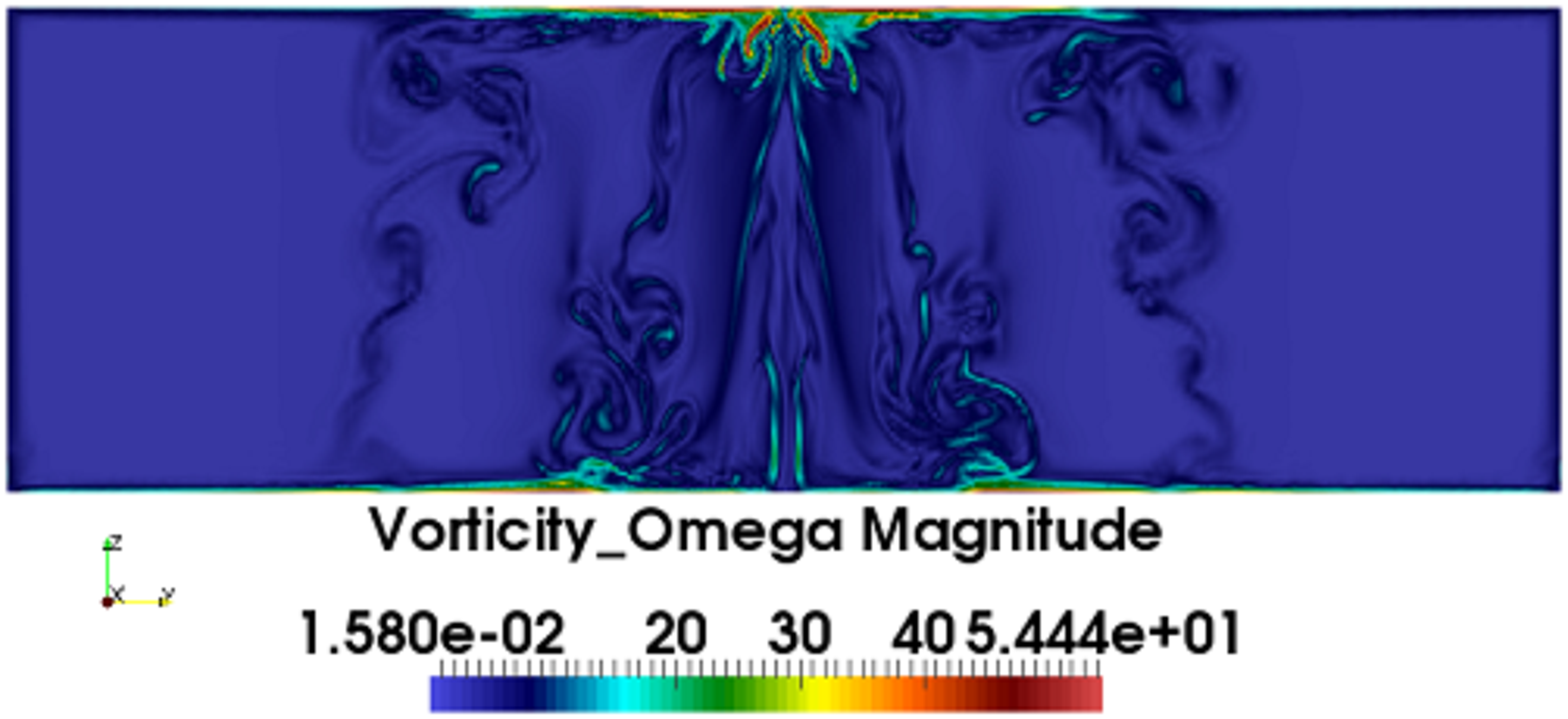}
\caption{Contours on the plane $x_1=0$ of the velocity magnitude
$|v|$ (first for $t=0.4$, second for $t=1.4$) and the vorticity magnitude $|\omega|$ (third for $t=0.4$, fourth for $t=1.4$) in the swirl case with $Re = 50,000$.}
\end{center}
\end{figure}

\section{A mathematical support of ``increasing velocity near the axis of symmetry and the lower boundary" }

When someone observes unstabilizing effects (such as increasing velocity in the time variation), it had better to compare with the stabilizing effects (such as regularity results mensioned in the introduction).  In the numerical computation in the previous section, Giga-Hsu-Maekawa's criterion (G-H-M's criterion in the following) was in mind. G-H-M's criterion (see (\ref{ghm}) also \cite[Theorem 1.3.]{GHM}) says that there are no type I blow-up solutions to (\ref{NS}) under a continuous alignment condition on the vorticity direction in the half space with no-slip boundary condition. More precisely, if the vorticity direction is controlled in some sense at the places where the absolute value of the vorticity is large, the flow is regular. Conversely, if we can construct solutions which do not satisfy such kind of continuous conditions at those places, we have more chances to find out the clue for possible blow-up solutions. Even in this paper we could not construct a blow up solution, this kind of thought might help us to explain or to predict some phenomenons with unstabilizing effects (such as increasing velocity phenomenons in this paper). Note that in our numerical computation, even we use the cylindrical domain instead of the half space, the behavior near the center of the lower (and upper) boundary is expected to be similar to the boundary of the half space under the same boundary condition.
Now we give a mathematical support why the maximum value is increasing near the saddle point on the boundary.
Let $\bar u=u_\theta e_\theta+u_re_r$ and $\bar u^{\perp}=u_re_\theta-u_\theta e_r$. We can calculate the vorticity on the boundary as  $\omega=\partial_z\bar u^{\perp}$.
Note that $\partial_zu_r(r,0)\to 0$ as $r\to 0$ due to the symmetry.
 We  can show that the direction of the vorticity $\omega/|\omega|$ is not continuous at the saddle point. In this case we need to assume $\partial_zu_\theta\not=0$ (on the boundary) and $\partial_ru_\theta\not=0$ (on the axis) near the saddle point.
These nonzero conditions may express ``shear flow effect by the swirl".
The vorticity along the axis is expressed as $\omega=2\partial_ru_\theta e_z$.
Thus $\omega/|\omega|$ along the $z$-axis and on the boundary is not continuous. Although, rigorously, it does not break the following continuous alignment condition in G-H-M's criterion:
\begin{equation}\label{ghm}
\left|\frac{\omega}{|\omega|}(t,x)-\frac{\omega}{|\omega|}(t,y)\right|\leq \rho(|x-y|),
\end{equation}
where $\rho(|x-y|)$ is any modulus continuous function, since the vorticity is zero at the saddle point. Nevertheless, by the numerical compuation, we see that there are high vorticity region near the saddle point and the boundary. We observe $|\omega|$ and three components of $\xi=\omega/|\omega|$ under a line which contains a point $(0, 0.05, -0.125)$ on the boundary and parallels to $z$-axis at $t=0.4$, see figure 8 (and see figure 9 for the case that the line parallels to $x_2$-axis). The bottom-axis in figure 8 represents the distance from the lower boundary while the bottom-axis in figure 9 represents the distance from the saddle point. Note that the coordinate of the saddle point is $(0, 0, -0.125)$ and the axial velocity $u_z=0$ on the boundary.
At the point $(0, 0.05, -0.125)$ on the boundary, magnitude of the vorticity $|\omega|$ is 60, and it clearly attains maximum value along the line in figure 8.
Moreover $\omega/|\omega|$ is also drastically changing (highly oscillating) near the boundary (see figure 8 and figure 9).
It shows some kind of unstability of $\xi$ at the location of the large vorticity magnitude.

\begin{figure}\label{xi}
\includegraphics[width=7.0cm
,keepaspectratio
]{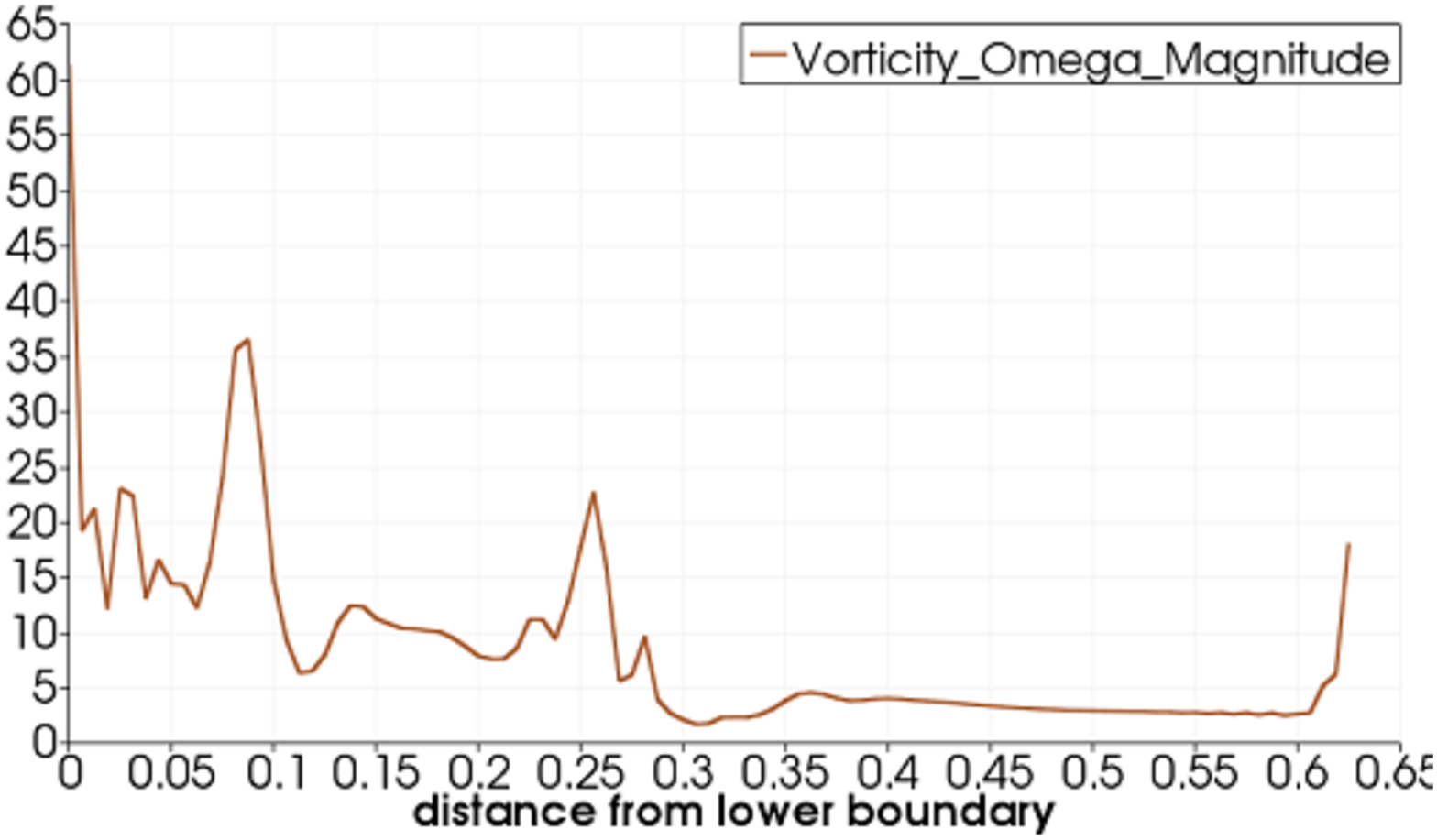}
\includegraphics[width=7.0cm
,keepaspectratio
]{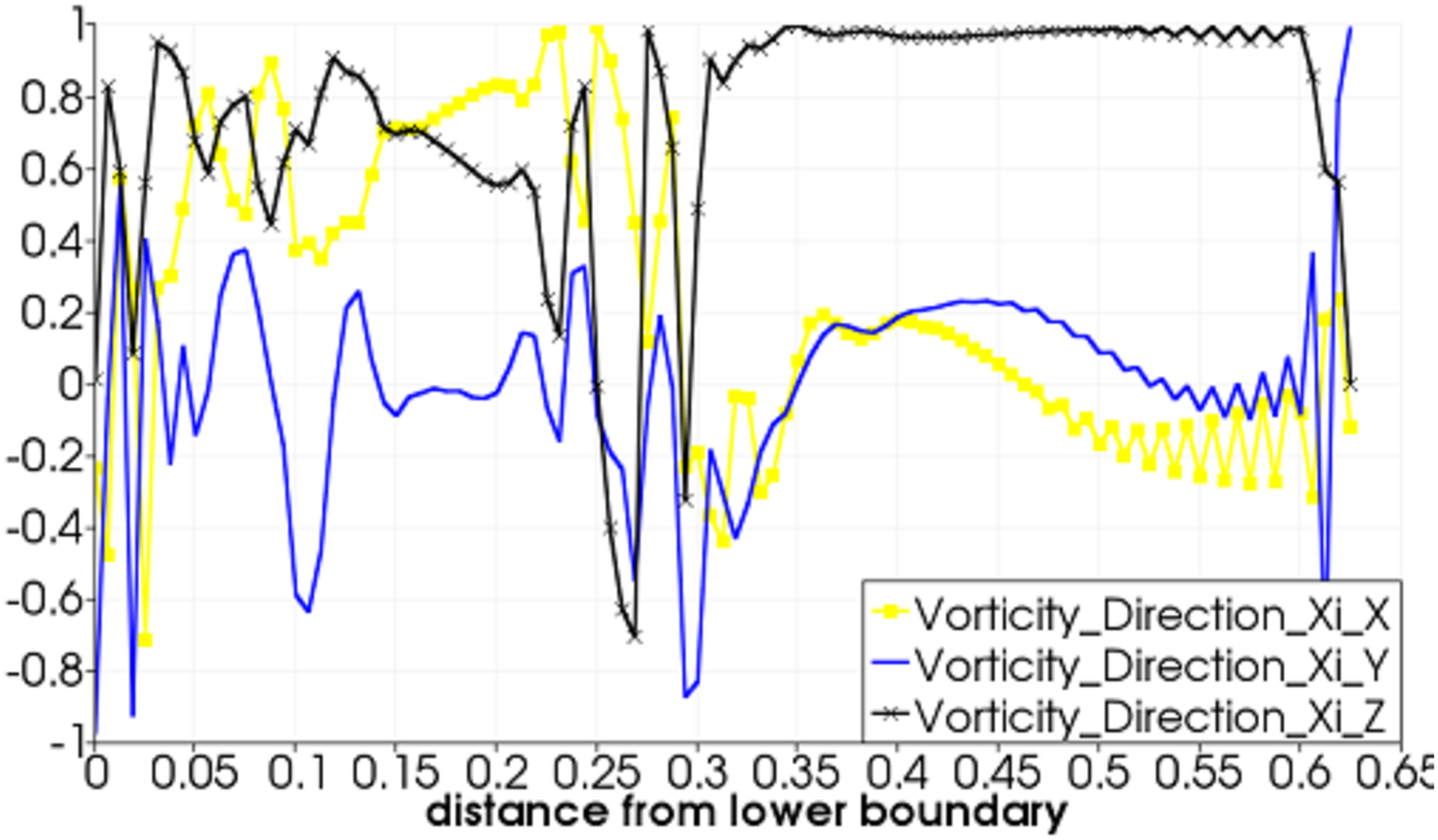}
\caption{Graphs of the quantities $|\omega|$ (top) and the three
components of $\xi$ (bottom) at $t = 0.4$ on the line parallel to the
$z$-axis through the point $(0, 0.05, -0.125)$.}
\end{figure}

\begin{figure}\label{xiy}
\begin{center}
\includegraphics[width=7.0cm
,keepaspectratio
]{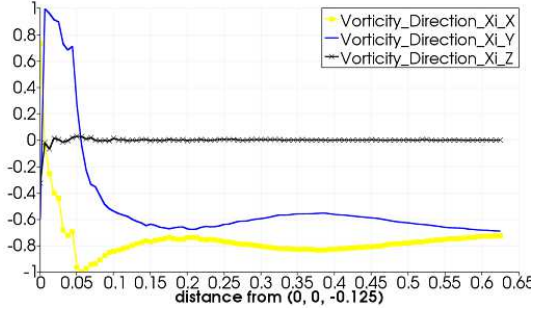}
\caption{Graphs of the three
components of $\xi$ at $t = 0.4$ on the line parallel to the
$x_2$-axis through the point $(0, 0.05, -0.125)$.}
\end{center}
\end{figure}


\section{Conclusion}
From literatures, it is natural to consider that the hyperbolic flow with swirl and saddle point on the boundary might be the key structure of flow and probable place for unstability effects occur near the no-slip flat boundary.
We showed the clear structure for the axisymmetric hyperbolic flow with swirl and observed the following phenomenons which are distinctly different from those without swirl:
(1) The distance between the maximum point of the velocity and the z-axis is drastically changing around some time we called it turning point.
(2) The velocity increases and obtains its extreme value (maximum) near the axis of symmetry and the boundary when time is close to the turning points.
By comparing these with the studies of tornadoes, it might help our understanding about the behavior of the velocity of wind near the ground which is very significant in the researches of tornadoes for reducing the damage cause by tornadoes or similar phenomenons.
(3) The downward flow near $z$-axis is observed. (See figure 6 the time variation diagram.) The downward wind inside the core of real tornado is also observed in two-celled vortex structure of the studies of numerical simulations for time-averaged velocity, see \cite[figure 4(b)]{IOT}. By comparing our observation with studies of tornadoes might enhance our understanding about the behavior inside the core of tornado for high swirl ratio. Furthermore, we also observed the properties of vorticity and its directions (see Section 3) which seems few in the literature of the numerical simulations of tornadoes.
Those phenomenons might be the clues for helping our understanding on the behavior near the saddle point and near the boundary. And our numerical approach might be helpful for other boundary shapes and other equations related to fluid mechanics in the future work. It might be useful for the studies of tornadoes arise or pass by a different landforms instead of flat plane.

\noindent
{\bf Acknowledgments.}\

At first, we would like to thank the anonymous referees for careful reading and valuable comments.
We would like to send our great gratitude to Professor Norikazu Saito for giving us valuable suggestions on the construction of the initial data.
We are also grateful to Professor Vladimir Sverak for letting us know interesting literatures \citep{CHKLSY, K}.
P.-Y.H. and T.Y. are supported by the "Program to Promote the Tenure Track System" of the Ministry of Education, Culture, Sports, Science and Technology.
Moreover,
P.-Y.H. gratefully acknowledges the support by the Iwanami Fujukai Foundation.
H.N. is supported by JSPS KAKENHI Grant Numbers 26800091 and 24224004,
by JSPS the Japanese-German Graduate Externship (Mathematical Fluid
Dynamics) and by Waseda University Project research of Research
Institute for Science and Engineering.
T.Y. is partially supported by JSPS KAKENHI Grant Number 25870004.
The numerical computation is based on Mr. Hiroaki Ishida's Master thesis at Hokkaido University on 2014 March. This research was supported by JST, CREST.

\appendix
\section{}\label{appA}
In this section, we present additional numerical results in order to
show two things.
The one is that the dependency of the numerical results in Section 2
on the discretization parameters $h$ and $\tau$ is qualitatively
small, and the other is the reason of the choice of the range of $z$
in $\Omega$, i.e., $-a<z<4a$. 

Firstly, we consider the former thing. We compute the no swirl and swirl cases for $Re = 50,000$, $10,000$
and $5,000$ by the stabilized Lagrange-Galerkin scheme with a coarse
mesh and a large time increment.
The maximum, minimum and average mesh sizes of the coarse mesh are
$2.48 \times 10^{-2}$, $1.93 \times 10^{-3}$ and $1.20 \times
10^{-2}$, respectively, where the strategy of the mesh generation is
the same, i.e., the mesh size around $z$-axis is smaller than that of
other part.
The time increment is set as $\tau = 1.66 \times 10^{-2}$.
In the following we call the numerical results in Section 2 ``results
A'' and the numerical results by the coarse mesh and the large time
increment ``results B''.
We compare results B with results A for $Re = 50,000$, $10,000$ and $5,000$.
Figure 10 shows graphs of maximum values of $|v|$ versus $t$ (left) and
the distance from the maximum point of $|v|$ to the $z$-axis versus
$t$ (right) in the swirl case for $Re=50,000$ (top), $10,000$ (middle)
and $5,000$ (bottom), where red and green colors are employed for
results A and results B, respectively.
Figure 11 displays corresponding graphs in the no swirl case for $Re=50,000$.
From figures 10 and 11, we can see that the two graphs in each figure
are qualitatively similar, while there is difference quantitatively.
Here, the graphs for $Re=10,000$ and $5,000$ in the no swirl case are
omitted, since they are also similar qualitatively.
We display additional information of results B in figure 12, which
shows cross-sections of $|v|$ (top) and $|\omega|$ (bottom) at $t =
0.4$ (left) and $1.4$ (right), for $Re = 50,000$ in the swirl case.
By comparing figure 7 (results A) with figure 12 (results B) we can
see that the behaviors are almost the same, although the magnitudes
are slightly different.
These results imply the former thing, i.e., the
dependency of the numerical results in Section 2 on the discretization
parameters $h$ and $\tau$ is qualitatively small.
\begin{figure}\label{meshAB}
\includegraphics[width=7.0cm
,keepaspectratio
]{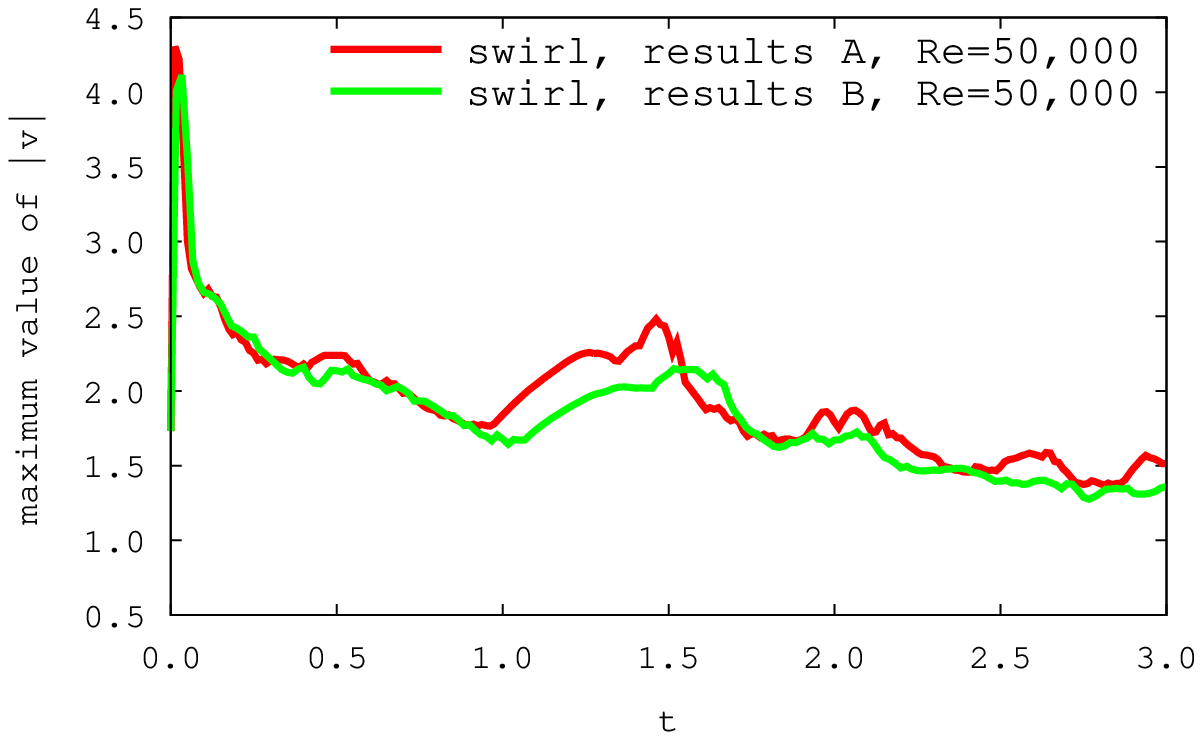}
\includegraphics[width=7.0cm
,keepaspectratio
]{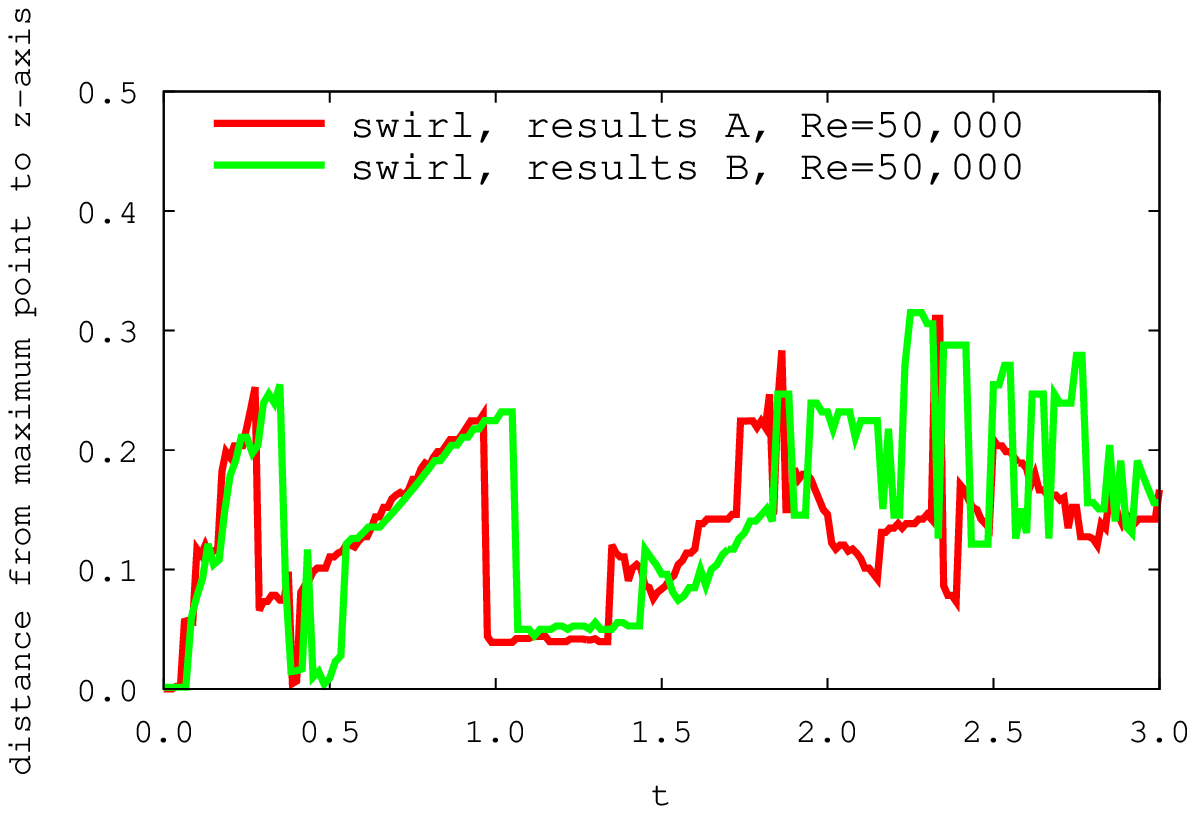}
\includegraphics[width=7.0cm
,keepaspectratio
]{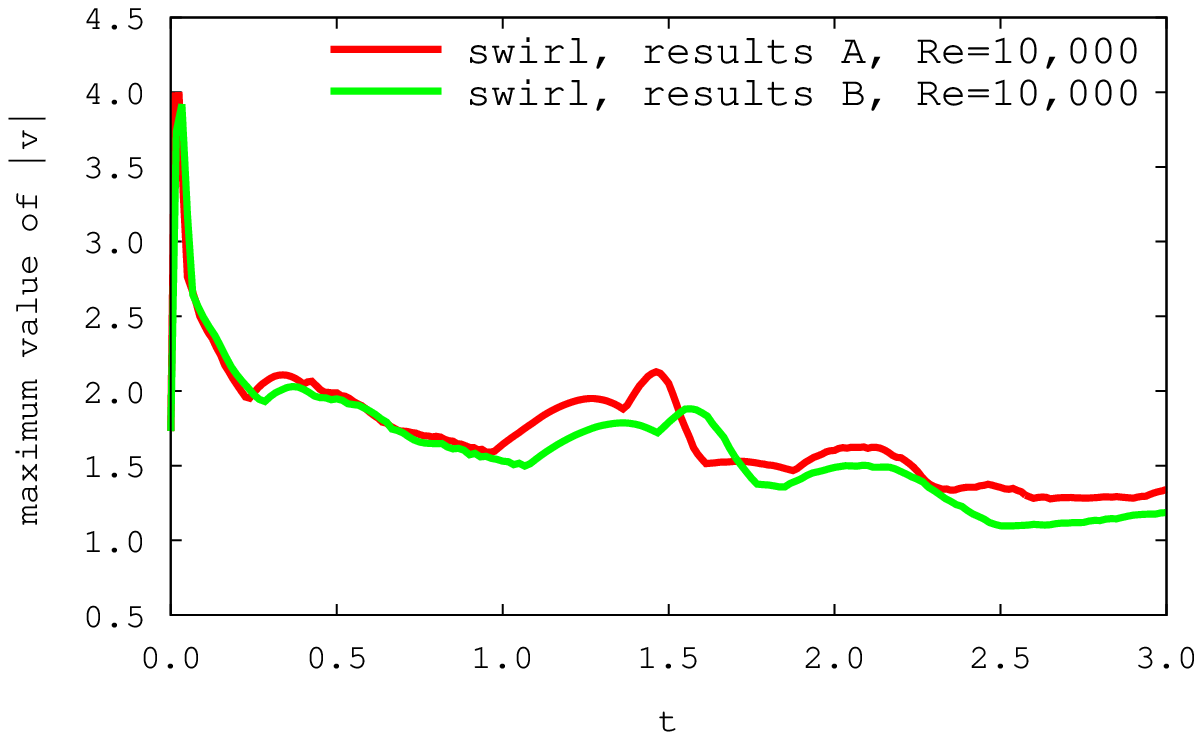}
\includegraphics[width=7.0cm
,keepaspectratio
]{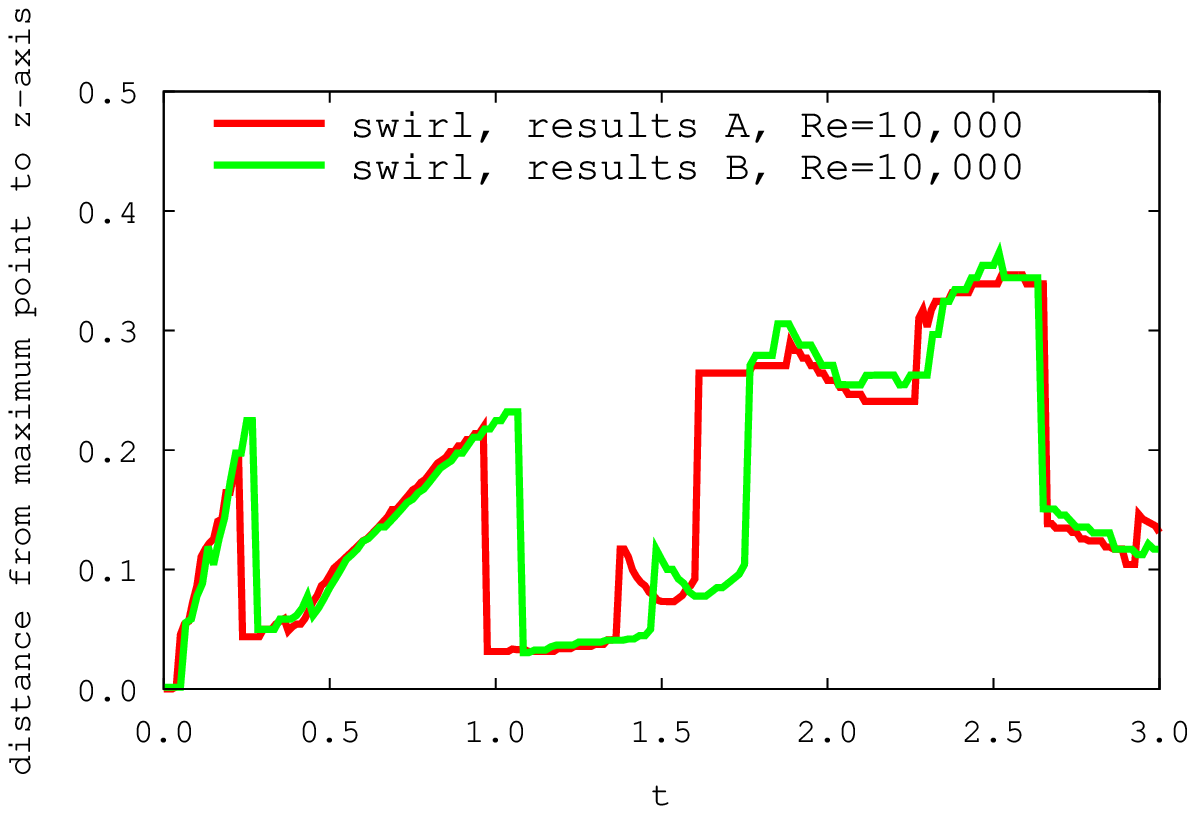}
\includegraphics[width=7.0cm
,keepaspectratio
]{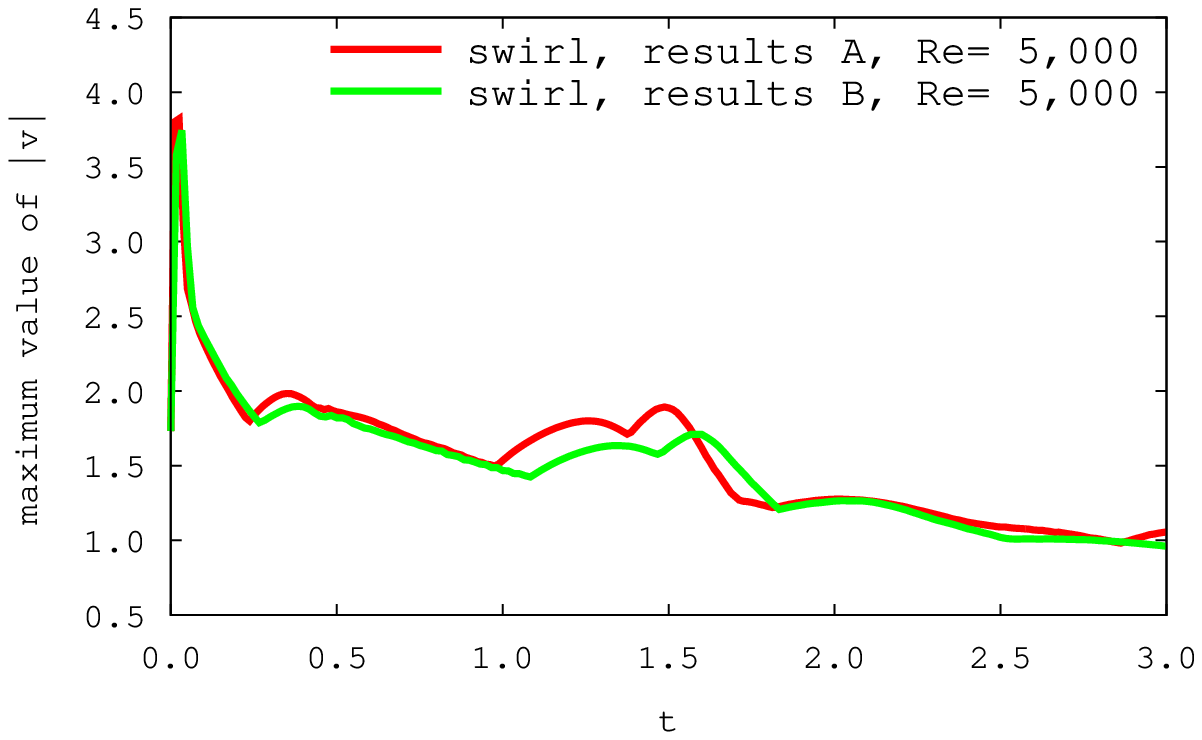}
\includegraphics[width=7.0cm
,keepaspectratio
]{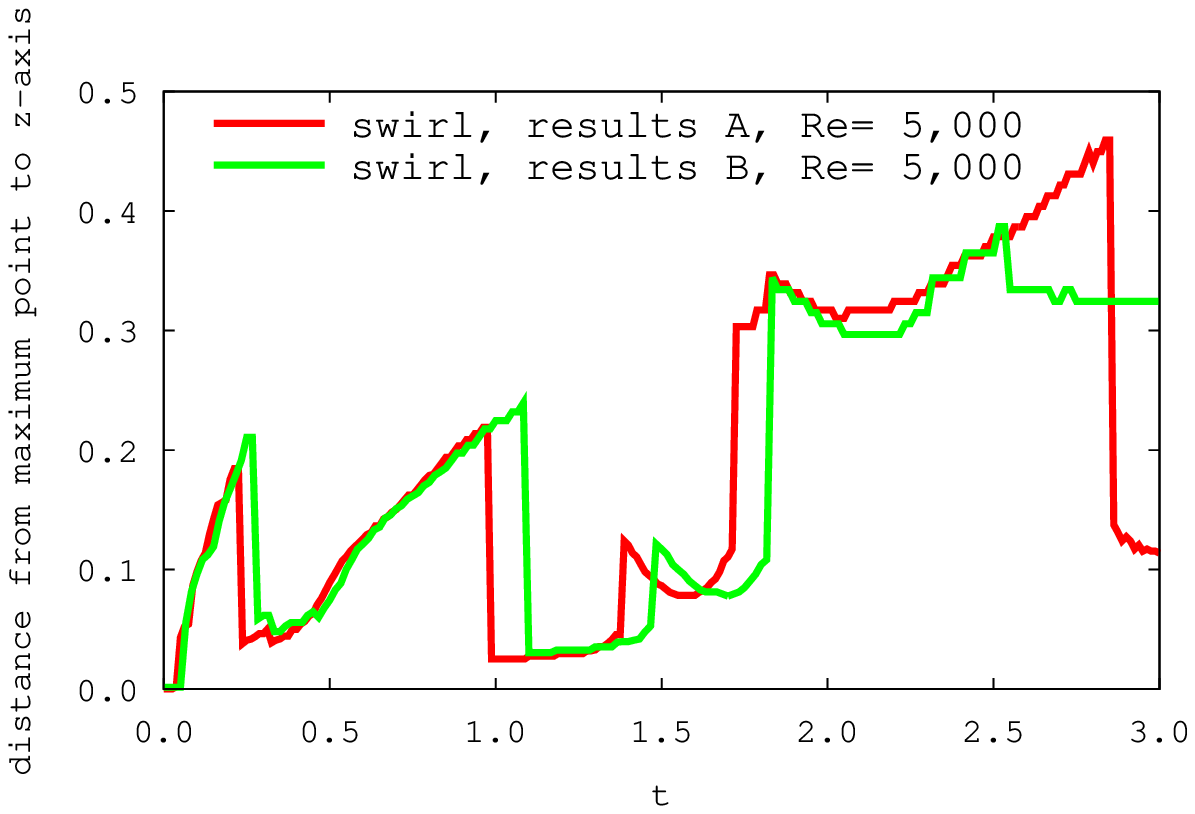}
\vspace{0.5cm}
\caption{Graphs of maximum values of $|v|$ versus $t$ (left) and
the distance from the maximum point of $|v|$ to the $z$-axis versus
$t$ (right) in the swirl case for $Re=50,000$ (top), $10,000$ (middle)
and $5,000$ (bottom), where red and green colors are employed for
results A and results B, respectively.}
\end{figure}

\begin{figure}\label{meshABnoswirl}
\includegraphics[width=7.0cm
,keepaspectratio
]{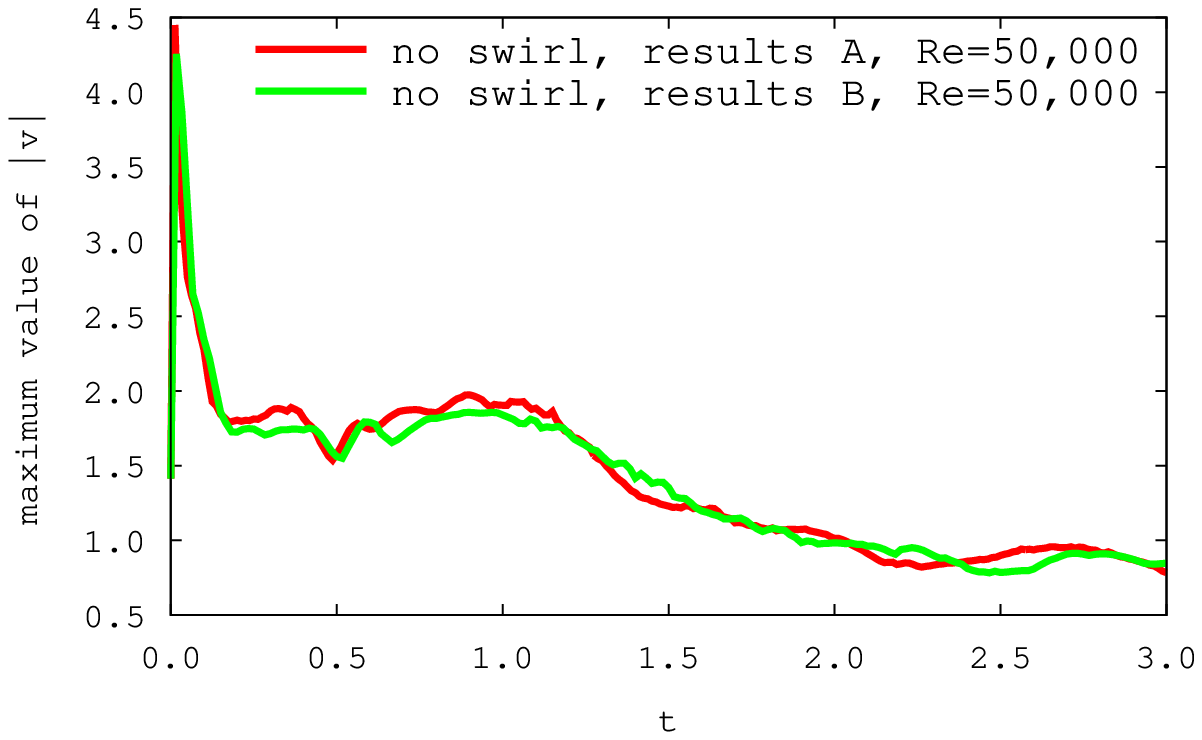}
\includegraphics[width=7.0cm
,keepaspectratio
]{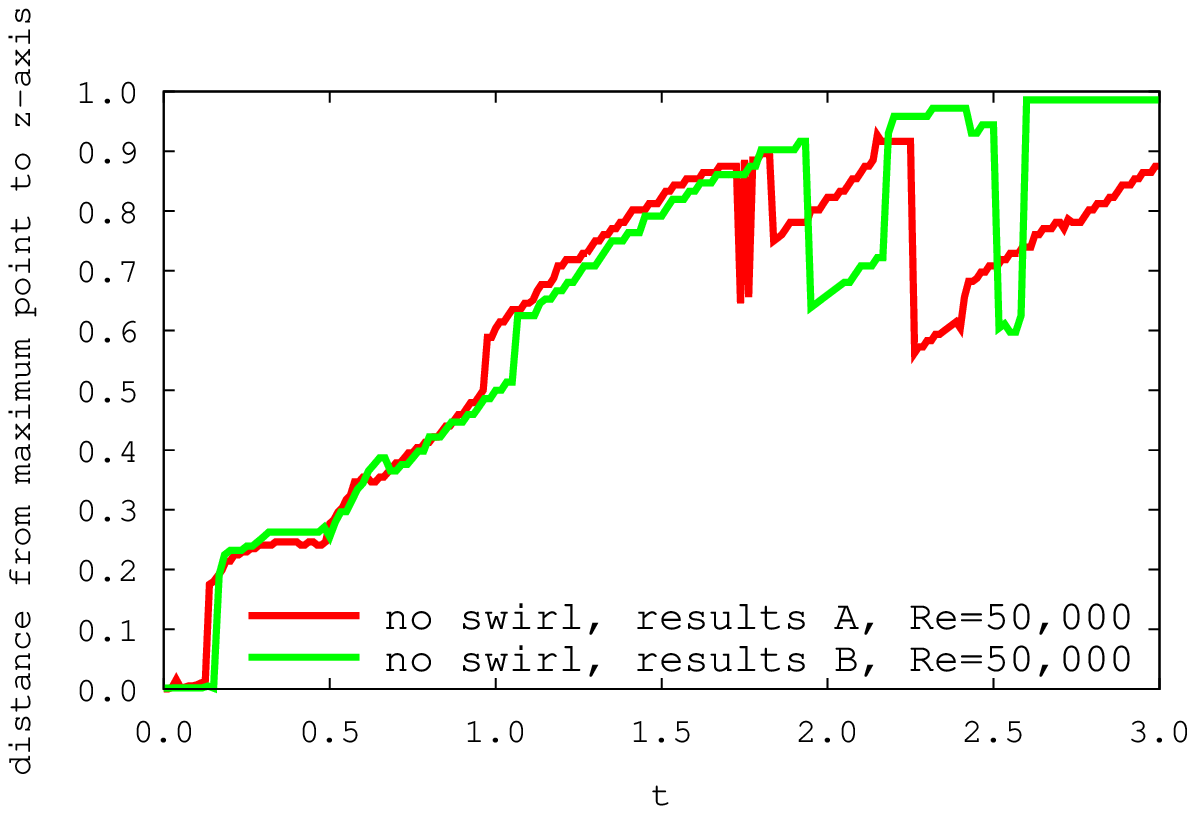}
\vspace{0.5cm}
\caption{Graphs of maximum values of $|v|$ versus $t$ (left) and
the distance from the maximum point of $|v|$ to the $z$-axis versus
$t$ (right) in the no swirl case for $Re=50,000$, where red and green colors are employed for
results A and results B, respectively.}
\end{figure}

\begin{figure}\label{keytimemeshB}
\includegraphics[width=7.0cm
,keepaspectratio
]{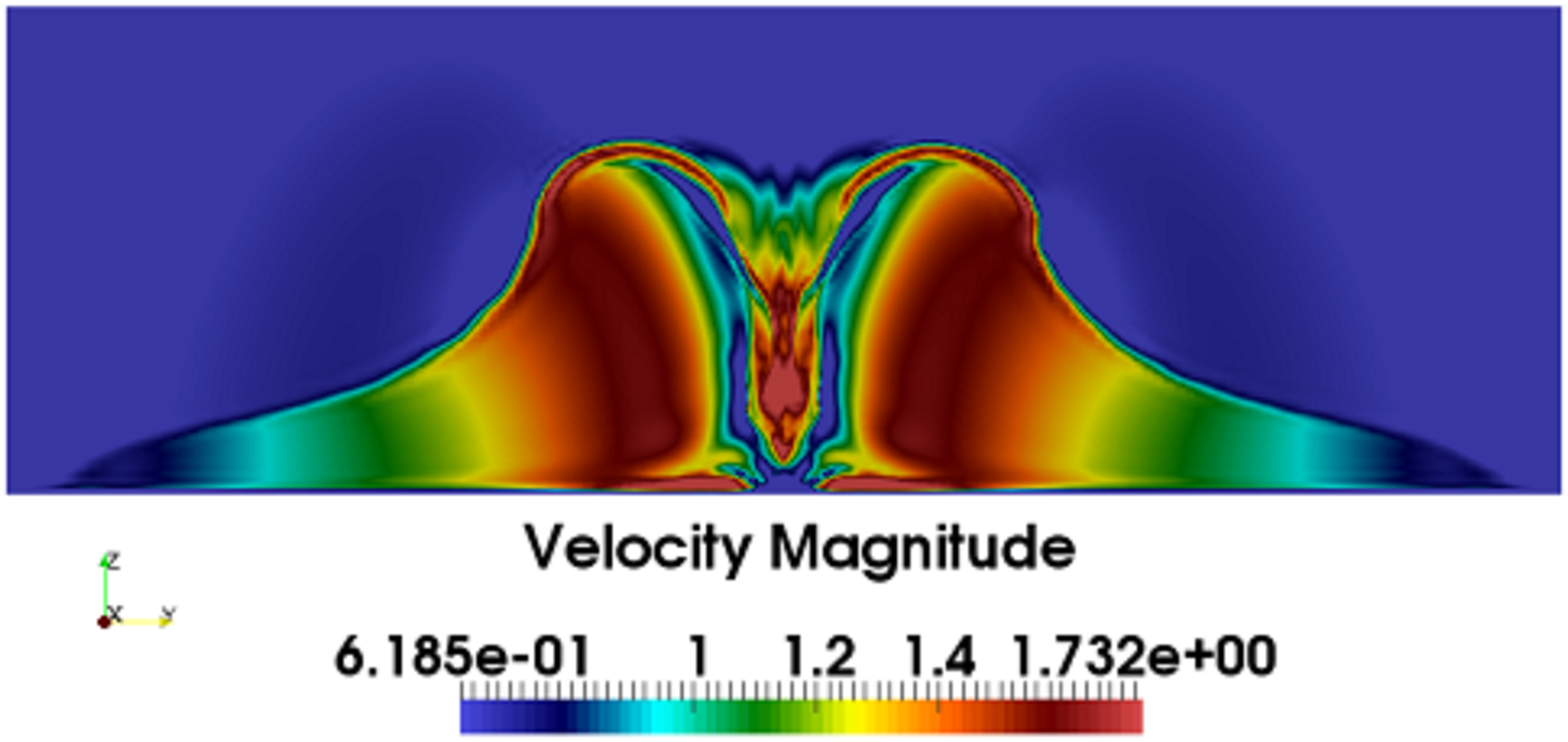}
\includegraphics[width=7.0cm
,keepaspectratio
]{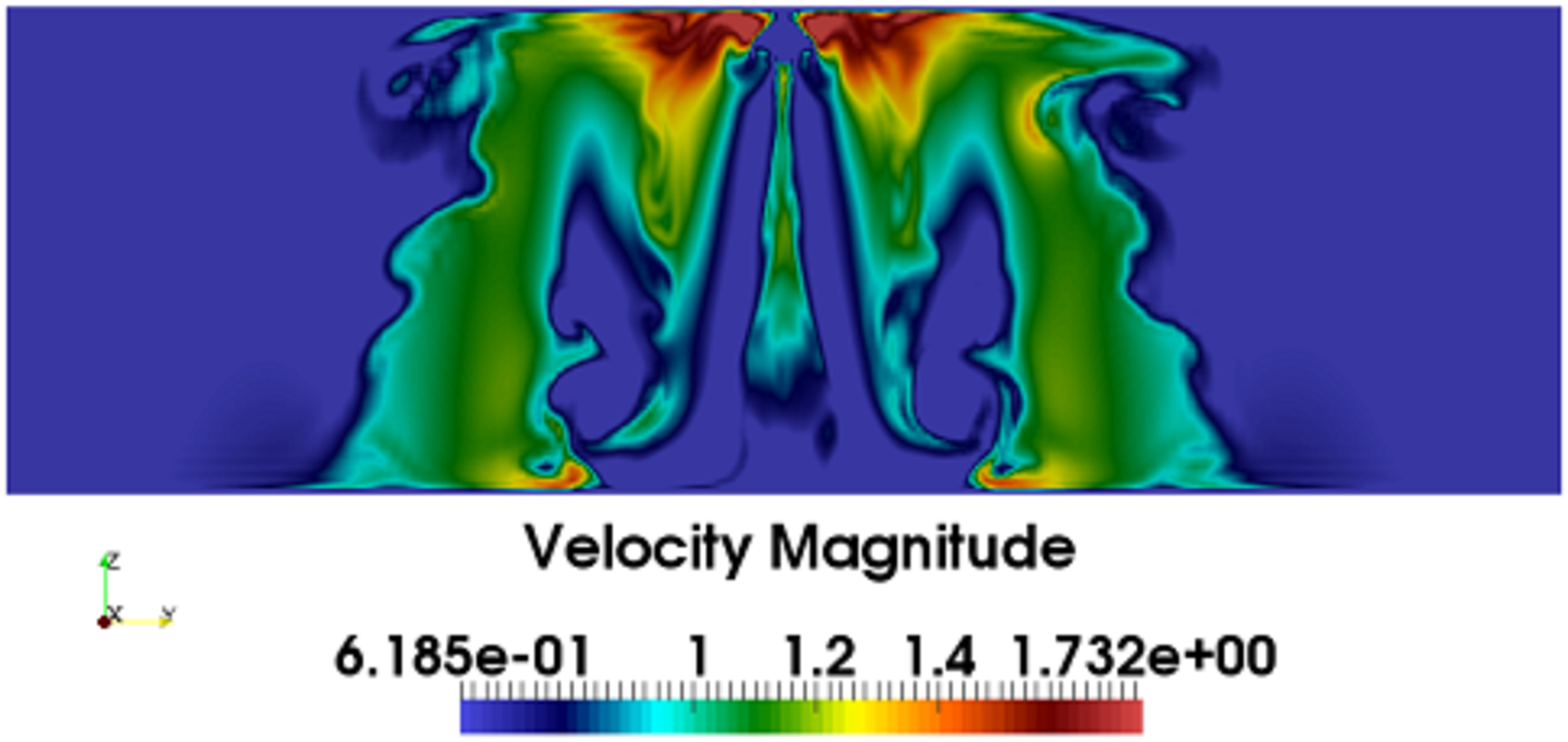}

\vspace{0.2cm}
\includegraphics[width=7.0cm
,keepaspectratio
]{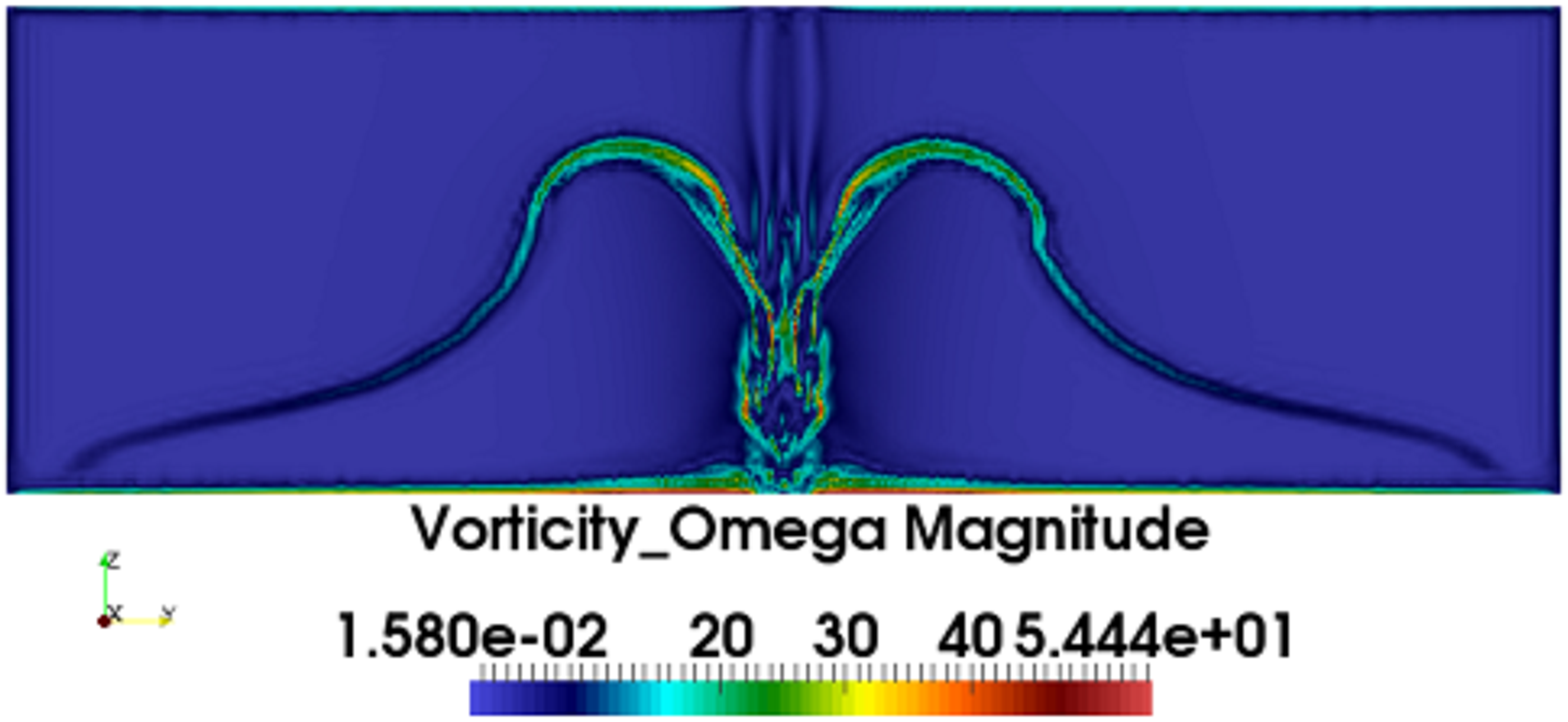}
\includegraphics[width=7.0cm
,keepaspectratio
]{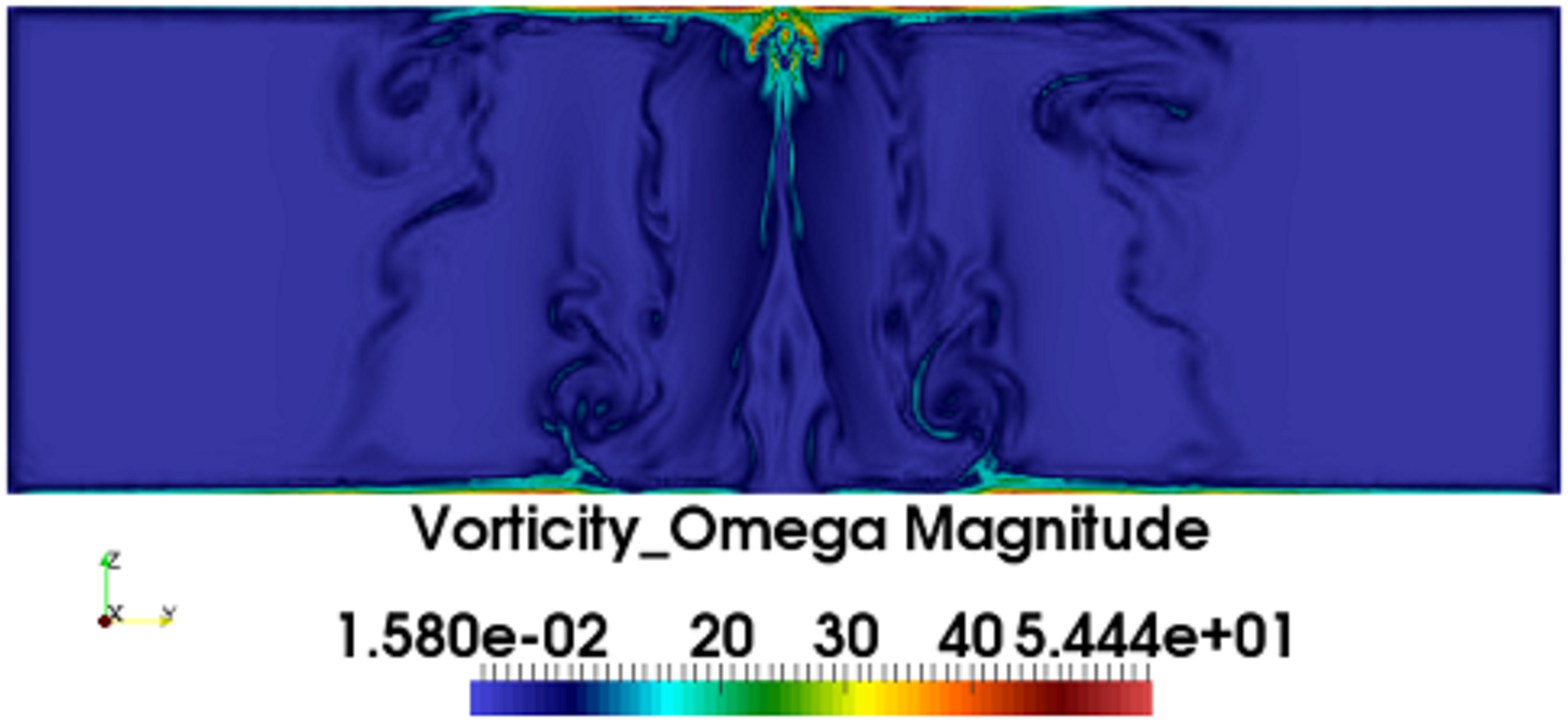}

\vspace{0.5cm}
\caption{Contours on the plane $x_1=0$ of the velocity magnitude
$|v|$ (top) and the vorticity magnitude $|\omega|$ (bottom) in the
swirl case with $Re = 50,000$ for results B. $t = 0.4$ (left) and $1.4$ (right).}
\end{figure}

Now we consider the latter thing.
Let $\Omega^\prime$ be a domain defined by $\Omega^\prime := \{ (x_1,
x_2, z) \in \mathbb{R}^3: -5a/2 < z < 5a/2,\ \sqrt{x_1^2+x_2^2} < 1
\}$.
We note that centers of $\Omega$ and $\Omega^\prime$ are $(0, 0,
3a/2)$ and the origin, respectively.
We compute the swirl case in the domain $\Omega^\prime$ for $Re =
50,000$ by the stabilized Lagrange-Galerkin scheme with a coarse mesh
and a large time increment, where the initial data is set by the same
functions in \eqref{initial}-\eqref{initialend}.
The mesh sizes and the time increment are the same as used in results B.
(The mesh is generated by a translation of the mesh for results B.)
We call the numerical results ``results C''.
We compare results C with results B in figure 13, which shows graphs
of maximum values of $|v|$ versus $t$ (left) and the distance from the
maximum point of $|v|$ to the $z$-axis versus $t$ (right) in the swirl
case for $Re=50,000$, where red and green colors are employed for
results B and results C, respectively.
It is observed that the ``increasing velocity phenomenon'' in results
B is clearer than that in results C, and that drastic changes of the
distance from the maximum point of $|v|$ to the $z$-axis are observed
in both results.
Hence, $\Omega$ is better than $\Omega^\prime$ in order to see both of
the ``increasing velocity phenomenon'' and the ``drastic changes''.
These results explain the latter thing, i.e.,
the reason of the choice of the range of $z$ in $\Omega$, $-a<z<4a$.

\begin{figure}\label{twoinitial}
\includegraphics[width=7.0cm
,keepaspectratio
]{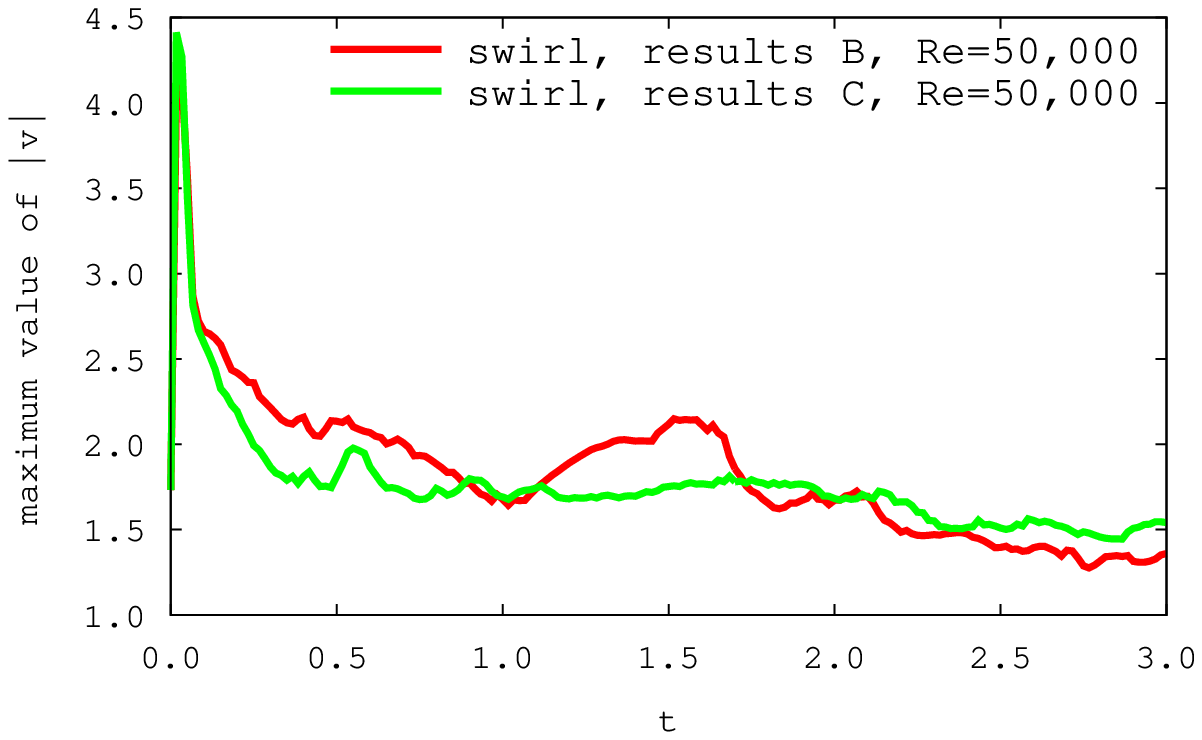}
\includegraphics[width=7.0cm
,keepaspectratio
]{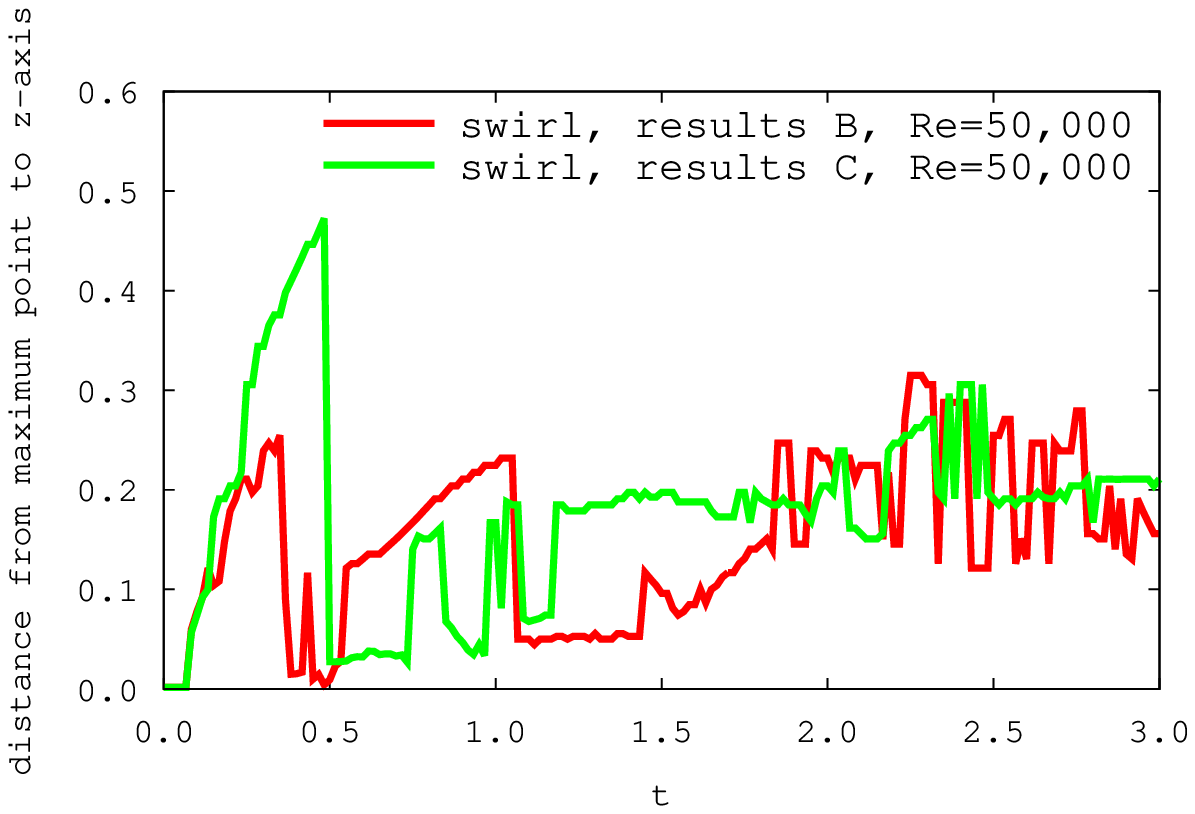}
\caption{Graphs
of maximum values of $|v|$ versus $t$ (left) and the distance from the
maximum point of $|v|$ to the $z$-axis versus $t$ (right) in the swirl
case for $Re=50,000$, where red and green colors are employed for
results B and results C, respectively.}
\end{figure}


\bibliographystyle{jfm}
\bibliography{HNY}

\end{document}